\documentclass[12pt]{article}

\usepackage{amsmath}
\usepackage{amssymb}
\usepackage{amscd}
\usepackage{graphicx}
\usepackage{psfrag}
\usepackage{bbm}

\usepackage{latexsym}
\usepackage{verbatim}
\usepackage{epsfig}
\usepackage{changebar}
\usepackage{graphicx}

\newtheorem{theorem}{Theorem}[section]
\newtheorem{lemma}[theorem]{Lemma}

\newtheorem{corollary}[theorem]{Corollary}

\newenvironment{proof}[1][Proof]{\begin{trivlist}
\item[\hskip \labelsep {\bfseries #1}]}{\end{trivlist}}
\newenvironment{definition}[1][Definition]{\begin{trivlist}
\item[\hskip \labelsep {\bfseries #1}]}{\end{trivlist}}

\newenvironment{remark}[1][Remark]{\begin{trivlist}
\item[\hskip \labelsep {\bfseries #1}]}{\end{trivlist}}

\newcommand{\qed}{\nobreak \ifvmode \relax \else
      \ifdim\lastskip<1.5em \hskip-\lastskip
      \hskip1.5em plus0em minus0.5em \fi \nobreak
      \vrule height0.75em width0.5em depth0.25em\fi}

\def\C{\mathbb  C}
\def\R{\mathbb
 R}

\def\P{\mathbf
 P}
\def\B{\mathbf
 B}

\def\S{\mathbf
 S}
\def\D{\mathbf
 D}

\def\e{\varepsilon}
\def\a{\alpha}
\def\mathbb{\mathbf
}

\def\a{\alpha}

\def\t{\theta}

\def\l{\ell}
\def\L{\mathcal L}
\def\nn {\vskip 0.3cm \noindent }
\def\n {\noindent}
\def\v {\vskip.2cm}
\begin{document}

\title{ Fibred Multilinks and singularities { $f \overline g$}}
\author{Anne Pichon and Jos\'e Seade}

\maketitle
\vskip0,2cm\noindent

%--------------------------
Anne Pichon, Institut de Math\'ematiques de Luminy,
UMR 6206 CNRS, Campus
de Luminy - Case 907, 13288 Marseille Cedex 9, France
\par\noindent
pichon@iml.univ-mrs.fr
\vskip0,2cm                                                                                                                                                                                                                                                                                                                                     
Jos\'e Seade, Instituto de Matem\'aticas, Unidad Cuernavaca,
Universidad Nacional Aut\'onoma de M\'exico, A. P. 273-3,
Cuernavaca, Morelos, M\'exico.
\par\noindent
 jseade@matcuer.unam.mx

\vskip0,3cm\noindent
Research partially supported by CNRS (France), CONACYT (Mexico) grant
G36357-E, a CONACYT-CNRS (France) Cooperation agre\-ement, 
and DGPA-UNAM (Mexico).
\vskip0,3cm\noindent

\vskip0,3cm\noindent

%-------------------------
 
\section{Introduction}
 
 Milnor's fibration theorem   in \cite {Mi1, Mi2} states that if 
$(\R^{n+k},0)  \buildrel{f}\over{\to} (\R^k,0)$ is the germ of a real analytic map with an isolated critical point at the origin, then for every sufficiently small sphere $\S_\e = \partial \B_\e$ around $0 \in \R^{n+k}$ one has that
the complement $\S_\e \setminus K$ of the link $K = f^{-1}(0) \cap \S_\e$ fibres over the circle $\S^1$. 
The proof of this result is by noticing first that for $\delta >0$ sufficiently small the tube $f^{-1}(\S^1_\delta) \cap \B_\e$ fibres over the circle $\S^1_\delta \subset \C$ of radius $\delta$, and then constructing a vector field that "inflates" this tube taking it into the complement of (a regular neighbourhood of) the link
in the sphere. When the map $f$ is from $\C^n$ into $\C$ and is holomorphic, Milnor shows that one actually has a much richer structure:

\v \n {\bf i)}  first, one does not actually need to have an isolated critical point of
$f$ to have such a fibration ;

\v \n {\bf ii)} in this case the projection map $\phi: \S_\e \setminus K \to \S^1$ can be taken to be
the obvious map $\phi = f/ \vert f  \vert$  (as Milnor shows in his book, this statement is false in general when $f$ is not holomorphic, even if one does have a fibration);

\v \n {\bf iii)}  the fibres $F$ of $\phi$ are parallelizable manifolds with the homotopy type of a CW-complex
of dimension $(n-1)$, and the link $K$ is $(n-3)$-connected.

\v \v
 
The geometry  of these fibrations associated to holomorphic  singularity germs has given rise to a rich literature,  as for instance   the  
theory of fibred knots and links,  open book decompositions, the results of Lawson and others about codimension 1 foliations, etc.
 This is not the case for the real analytic germs.  There are various reasons for this, in particular because it is difficult to find examples of real analytic singularities with an isolated critical point, and it is even harder to study their underlying geometry (the topology of the link and of the fibres, the monodromy,...).   Several natural -related- problems arise:  i) find examples of real analytic germs with isolated singularities and describe their underlying geometry; ii) relax the conditions in Milnor's fibration theorems for real singularities in order  to include larger families; and  iii) study which fibred links (and multilinks) can be realized by real analytic singularities.

These questions have been addressed by several authors in various ways. For instance Looijenga in \cite {Lo} realized fibred links which do not come from complex singularities. In \cite{Pe} Perron realizes the eight knot as the link of a real analytic map $\R^4 \to \R^2$, and   Lee Rudolph \cite{Ru} gives another (simpler) algebraic construction of this knot. In \cite{AK} Akbulut and King  
proved that every knot (and every link) in a sphere can be realized by 
real algebraic equations; this result was later
improved in \cite{MVS}, getting a better description of the method for obtaining
 the specific equations for a given knot.
In \cite {AC} A'Campo showed that the real mapping
$\R^{4} \cong \C^2 \to \R^2 \cong \C$ given by $(z_1,z_2) \mapsto \bar z_1 \bar z_2 (z_1 + z_2) $ defines an isolated singularity with a Milnor fibration which is not equivalent to a fibration of a complex singularity. In \cite{Se} (see also \cite{RSV, Se2}) Seade proves that the singularities defined in $\R^{2n} \cong \C^n$ by
$z_1^{a_1} \bar z_{\sigma_1} + \cdots + z_n^{a_n} \bar z_{\sigma_n}$, $a_i >1$, where 
$(\sigma_1,\cdots,\sigma_n)$ is any permutation of $(1,\cdots,n)$, define isolated singularities 
with a Milnor fibration. In \cite{PS} is studied the case $h(z_1,z_2)  =  z_1^{a_1} \bar z_2 + z_2^{a_2} \bar z_1$
describing their topology and proving (together with \cite {Pi1}) that the corresponding fibration is equivalent to the fibration defined by 
$(z_1,z_2) \mapsto \bar z_1 \bar z_2 (z_1^{a_1+1}  + z_2^{a_2+1} ) $, which are singularities of the same type as A'Campo's example. 
This type of singularities also appear in Lee Rudolph's paper  \cite {Ru}, whose work suggested that if $f,g$ are holomorphic map-germs   $(\C^2,0) \to (\C,0)$ with an isolated critical point at $0$ and   no common branch, then the
 map $(\C^2,0) \buildrel{f \overline g}\over \to (\C,0)$ has an isolated critical point at $0$ if and only if  the link $L_{f \overline g} = L_f - L_g$ is a fibred link, where $L_f$ and $L_g$ are the the links of the corresponding singularities in $\C^2$, endowed with their natural orientations. 
This was proved by Pichon in \cite {Pi2}, where she also
 proved that if 
 $(\C^2,0) \buildrel{f \overline g}\over \to (\C,0)$ has an isolated critical point at $0$, then the underlying Milnor fibration  is an open book fibration of the oriented link $L_{f \overline g} $ given by the obvious map ${f \overline g}/\vert {f \overline g}\vert$ in a neighbourhood of $L_{f \overline g} $. This enables us to realize a large family of fibrations of plumbing links in $\S^3$ as the Milnor fibrations of real analytic germs ${f \overline g}$.

These results arose our interest in studying singularities of the type $f \overline g$ where $f, g$ are both holomorphic functions $(\C^n,0) \to (\C,0)$. An obvious problem is that if $n >2$ then the zero-sets of $f$ and $g$ will be hypersurfaces that necessarily meet in an analytic set of dimension $>0$; thus 
$0 \in \C^n$ may not possibly be an isolated critical point of  $f \overline g$. However simple examples
(see section 1 below) show that $0 \in \C$ is often an isolated critical value of $f \overline g$.  This motivated the following theorem, which says that Milnor's fibration theorem for complex singularities,
 together with statements (i)-(iii) above, remains valid for these singularities:

\nn {\bf Theorem 1.} {\it Let $X$ be a complex analytic variety of dimension $n > 0$ in an open set in $\C^N$, with an isolated singularity at $0$, and let
$f,g : (X,0) \to (\C,0)$ be holomorphic maps such that the real analytic map
$$f \overline g: (X,0) \to (\R^2,0)$$
has an isolated critical value at $0 \in \R^2$. Let $\L_X = X \cap \S_\e$ be the link of $0$ in $X$ (so it is a smooth manifold of real dimension $(2n-1)$), and let $L_{fg} = (fg)^{-1}(0) \cap \L_X$ be the link of
$fg$ in $\L_X$. 
Then the map:
\[\phi: = \frac{f \overline g}{\vert  f \overline g \vert} \,:\, \L_X \setminus L_{fg}  \longrightarrow \S^1 \subset \C\]
is the projection of a $\C^\infty$ (locally trivial) fibre bundle, whose fibres $F_\t$ are diffeomorphic to the complex manifolds $(f/g)^{-1}(t) \cap \buildrel{\circ} \over \D_\e$, where $t \in \C$ is a regular value of the meromorphic
function $f/g$ and $\buildrel{\circ} \over \D_\e$ is the interior of the disc in $\C^N$ whose boundary is $\S_\e$. Hence each fibre has the homotopy type of a CW-complex of dimension $(n-1)$. Moreover, if the tangent bundle of $X \setminus\{0\}$ is trivial (for instance if $X \cong \C^n$), then the fibres of $\phi$ are parallelizable manifolds.
}

\v
It is clear that
 the link $L_{f \overline g}$ of $f  \overline  g$ equals the link 
$L_{f g}$ of the complex singularity 
$\{fg =0\}$ as unoriented manifolds, hence $L_{f \overline g}$ is $(n-3)$-connected when $X$ is $\C^n$
(by \cite[5.2]{Mi2}). Notice also that on 
$\L_x \setminus L_{f \overline g}$ the map $\phi$ coincides with 
$\frac{f/g}{\vert f/g \vert} $, so Theorem 1 can be regarded as a fibration 
theorem for meromorphic map-germs.

\v \v

This theorem is also valid and interesting when $n=2$, and  the investigation of that situation takes sections 4-7 below. We show that the statement of Theorem 1 can be made much stronger in this case.  Notice that we do not ask here that $f$ and $g$ have an isolated critical point at $0$. This means that each branch of the curve 
$X \cap \{f^{-1}(0) \cup g^{-1}(0)\}$ may have a certain multiplicity $\ge 1$ and the intersection
$L \,=\, X \cap \{f^{-1}(0) \cup \overline g^{-1}(0)\} \cap \S_\e$ with a small sphere can be regarded as a multilink (c.f. \cite {EN}):
$$L \,= (n_1 K_1 \cup  \cdots \cup n_r K_r)\, \bigcup \,(- n_{r+1}  K_{r+1} \cup  \cdots \cup - n_{r+s} K_{r+s})\;,$$
where $K_1,\cdots,K_r$ are the components corresponding to the branches of $f$, 
$K_{r+1},\cdots, K_{r+s}$ those of $g$, the $n_i$ are the corresponding multiplicities, and for the components corresponding to $\overline g$ we take the negative multiplicity because they get the opposit orientation (see section 4 for details). Our main result in this direction is (Theorem 6.2):

\nn {\bf Theorem 2.} {\it 
Let $(X,p)$ be normal complex surface singularity and let $f: (X,p) \to ({\mathbf C},0)$ and $g: (X,p) \to ({\mathbf C},0)$ be two holomorphic germs without common branches. Then the following three conditions are equivalent:

\v \n {\bf (i)} the
real analytic germ $f \overline g : (X,p) \to (\mathbf R^2,0)$ has $0$ as an isolated critical value at $0$;

\v \n {\bf (ii)} the multilink $L_f -L_g$ is fibred; 

\v \n {\bf (iii)} the multiplicities corresponding to  $f$ and $g$ are distinct at each rupture vertex of the dual graph of a resolution of the curve $X \cap \{f^{-1}(0) \cup \overline g^{-1}(0)\}$.

\v
Moreover, if these conditions hold, then the Milnor fibration $\frac{f \overline{g} }{ |{f \overline {g}}|} : 
{\L}_X \setminus (L_f \cup  L_g)$ of ${f \overline g}$ is a fibration of the multilink $L_f - L_g$.
}

\v

The third condition above is very simple to check and gives an amazingly easy
criterion to decide whether or not the other two conditions are satisfied. 

The proof of Theorem 2 uses our Theorem 1 above, together with  
a generalization of a theorem of D. Eisenbud and W. Neumann \cite {EN} concerning fibred multilinks in general (section 4),  the study in section 5 of the geometry of the map $\frac{f \overline{g} }{ |{f \overline {g}}|}$ near the link $L_{f \overline g}$, similar to that in \cite {Pi2}, and the use (in section 6) of the 
discriminantal ratios introduced in \cite {LMW}.

In section 7 we use the previous results and the main theorem in \cite {Pi3} to prove the following realization theorem for multilinks:

\nn {\bf Theorem 3.} {\it  Let $M$ be a plumbing  $3$-manifold  and let $L_1$ and $L_2$ be two plumbing multilinks in $M$ with positive multiplicities which verify the following conditions : 
 
1.  the three multilinks $L_1 \cup (- L_2)$, $L_1$ and $L_2$ are fibred;

2. $M$ is either $\S^3$ or 
 homemorphic to the link ${\cal L}_X$ of a  taut  surface singularity $(X,p)$.
 \n 
Then, there exist two holomorphic germs $f,g:(X,p) \to ({\mathbf C},0)$ without common branches such that $L_1$ (and $L_2$) is the multilink of $f$ (respectively $g$); $L_1 \cup (-L_2)$ is the multilink associated to the real analytic germ $f \overline{g} : (X,p) \to ({\mathbf C},0)$; and
$$\frac{f \overline g}{\vert f \overline g \vert} \,:\, 
{\cal L}_X \setminus (L_1 \cup (- L_2)) \longrightarrow \S^1\;,$$
is a fibre bundle that realizes $L_1 \cup (- L_2)$ as a fibred multilink.
}

\v
We recall \cite{L2} that a complex surface singularity is taut when its topology determines the analytic structure. This
applies, for instance, to all rational double and triple points, as well
as to the cyclic (cusp) singularities of Hirzebruch.

This means in particular that $f \overline{g}$ has $0$ as an isolated critical value and that the Milnor fibration of $f \overline{g}$ is a fibration of the multilink  $L_1 \cup (-L_2)$. Notice that theorem 4.1 enables us to check  whether a plumbing multilink in a $3$-manifold is fibred or not. Moreover, when it is fibred, then the results of \cite{Pi1} enable us to  compute explicitly the genus of the fibre and the underlying monodromies. See for instance \cite{PS} for these computations in some examples.

The first part of this article, sections 1 to 3, are devoted to studying generalisations of Milnor's fibration theorem in all dimensions. Section 1 proves a fibration theorem for real analytic map-germs 
$(X,p) \to (\R^k,0)$ which have an isolated critical value at $0$ and 
satisfy Thom's $a_f$-condition for some Whitney stratification of the complex analytic space $X$;  we give several examples of such maps. This section is largely indebted to \cite{Mi1}, an unpublished article by J. Milnor which was pointing towards his fibration theorem in \cite{Mi2}  (finally proved in a different way).  In section 2 we prove that if the real analytic map is of the form
$f \overline g$ with $f, g$ holomorphic  complex valued functions and  $f \overline g$ has an isolated critical value at $0 \in \C$, then one has a fibre bundle:
$$ \phi = \frac{f \overline g}{\vert f \overline g\vert} : \L_X \setminus L_{fg} \longrightarrow \S^1$$
as in Theorem 1 above. In section 3 we prove that with the conditions of section 2 one always has a fibre bundle as in section 1, and this fibre bundle is equivalent to that of section 2. A consequence of this equivalence is that the fibres of $\phi$ are diffeomorphic to Stein manifolds in an affine space, and therefore the theorem of Andreotti-Frankel \cite {AF} implies the statements in Theorem 1 about their topology. The final statement about the parallelizability of the fibres follows from  general 
arguments given in \cite{KM}.

 A key point for the arguments in sections 2 and 3 is to notice that away from the link $L_{fg}$ the map
 $\phi = \frac{f \overline g}{\vert f \overline g\vert}$ equals the map 
 $\frac{f / g}{\vert f / g\vert}$ and therefore one can mimic the proof of Milnor in \cite{Mi2}, replacing  holomorphic functions by 
the meromorphic one $f/g$.

\section{ Milnor's fibration and the Thom $a_f$-condition}

Let $U$ be an open neighbourhood of $0$ in $\R^{n+k}$ and let 
 $X \subset U$ be a real analytic variety  with an isolated
singularity at  $0$. Let $\tilde f : (U,0) \to (\R^k,0)$ be a real analytic germ which is
 generically a submersion, {\it i.e.} its jacobian
 matrix $D \tilde f$ has rank $k$ on a
dense  open subset of $U$. We denote by $f$ the restriction of $\tilde f $ to $X$.
Let us recall some concepts that go back to the work of R. Thom and others (see for instance 
\cite {GM}). We say that $x \in X \setminus \{0\}$ is {\it a regular point} 
of $f$ if $Df_x:T_xX \to \R^k$ 
is a submersion. This means that $x$ is not a critical point of $\tilde f $ and 
the kernel of $D \tilde f_x$ is transverse to
the tangent space $T_xX$. A point in $X$ which is not a regular point of $f$ is said to
be {\it a critical point}. A point $y \in \R^k$ is {\it a regular value} of $f$ if there 
are no critical points of $f$ in $f^{-1}(y)$; otherwise we say that $y \in \R^k$ is 
{\it a critical value} of $f$, which means that there 
exists a critical point $x \in X$ such that $f(x) = y$. We say that $f$ has 
{\it an isolated critical value} at $0 \in \R^k$ if there is a neighbourhood $\D_{\delta}$ 
of $0$ in $\R^{k}$ so that all points $y \in \D_{\delta} \setminus \{0\}$ are regular values 
of $f$. In the last chapter of his book \cite {Mi2} Milnor proved a fibration theorem for real analytic singularities with an isolated critical point at the origin. Our aim in this section is to extend that result to (certain) real singularities with an isolated critical value at the origin. Let us motivate this with some examples.

\nn
{\bf 1.1 Examples. a)}  Consider first the map $f: \R^4 \cong \C^2 \longrightarrow  \R^2 \cong \C$  
defined by:
$$f(x,y) = \overline x \overline y (x^p + y^q)\;.$$
This can be written as $f = (f_1,f_2)$ where $f_1 = \frac{1}{2}(f + \overline f)$ and 
$f_2 = \frac{1}{2i}(f - \overline f)$ are its real and imaginary parts. A straight-forward
computation shows that $f$ has an isolated critical point at $0 \in \R^4$ if $p,q > 1$ and
not both of them are 2. In this case Milnor's condition to have 
a  fibration holds (c.f. \cite {PS, Pi2}).

\nn
{\bf b)} Consider now the map $\C^2 \buildrel{f}\over{\longrightarrow}  \C$  defined by
$f(x,y) = \overline x^3 (x^2 +y^3)$. In this case the  critical set of 
 $f$ is the axes $\{ x = 0\}$. Hence its  critical 
points are non-isolated, but $0\in \C$ is its only 
 critical value.  
 
\nn
{\bf c)}
Now let $h: \R^{2n} \cong \C^n \to  \R^2 \cong \C$, $n > 1$, be defined by $h = f \overline g$ where $f$ is the Pham-Brieskorn polynomial 
$$f(z_1,\cdots,z_n) = (z_1^{a_1} + \cdots z_n^{a_n})\,\,\,,\, a_i \geq 2, $$
and $g (z_1,\cdots,z_n) = \overline z_1 \cdots  \overline z_n$.
Writing   $h = (h_1,h_2)$, where $h_1$ and $h_2$ are its real and imaginary parts, a 
straight-forward computation shows that the rank of the corresponding jacobian matrix is less than 2 if and only if the following conditions are satisfied:
\v \n
{i) } $\quad \vert g \vert \,\vert a_i z_i^{a_i} \vert \,=\, \vert g \vert \, \vert f \vert\;, \quad i = 1, \cdots, n\;;$

\v \n
{ii) } $\quad a_i z_i^{a_i}  z_k \,=\, a_j z_j^{a_j}  z_k\;, \quad \hbox{for all distinct } i, j, k$ (if  $n =2$ just delete $z_k$ in both sides of the equation); and

\v \n
{iii) } $\quad a_i a_j z_i^{a_i}  \overline z_j^{a_j}  \vert g \vert^2 \,=\, \vert g \vert^2  \, \vert f \vert^2  \;$, for all distinct  $i, j$.
 \v \n For  points in $\C^n \setminus \{fg = 0\}$ the third set of equations follows from the others 
and the conditions reduce to:

\v \n
{i) } $\quad \vert a_i z_i^{a_i} \vert \,=\, \vert f \vert\;, \quad i = 1, \cdots, n\;;$ and 

\v \n
{ii) } $\quad a_i z_i^{a_i}   \,=\, a_j z_j^{a_j}  \;, \quad \hbox{for all distinct } i, j$.

 \v \n Then an easy computations shows that if the $a_i$ satisfy:
\[ \frac{1}{a_1 } + \cdots + \frac{1}{a_n } \, \ne \, 1\,,\]
then $0 \in  \R^2$ is an isolated critical value of $h$. Notice that for $n =2$ this condition is satisfied whenever one $a_i$ is more than 2. For n = 3 the condition is satisfied whenever 
the unordered triple $(p,q,r)$ is not $(2,3,6)$, $(2,4,4)$ or $(3,3,3)$.

\nn

 Let $\tilde f : (U,0) \to (\R^k,0)$ be again a real analytic map, and assume $f = \tilde f|_X$  has 
 an isolated critical value at $0 \in \R^k$. We set $V = f^{-1}(0) = \tilde f^{-1}(0) \cap X$.
According to \cite {Hi, LT}, there exists a Whitney stratification
$(V_\a)_{\a \in A}$ of $U$ adapted to $X$ and $V$. 

\nn
{\bf 1.2 Definition.} The  Whitney stratification $(V_\a)_{\a \in A}$ satisfies
{\it the Thom} $a_f$ {\it condition} with respect to the map $f$
 if for every sequence of points 
$\{x_n \} \in X \setminus V$ converging to a point $x_0$ in a stratum $V_\a \subset V$
 such that the sequence
of tangent spaces $T_{x_n}(f^{-1}(f(x_n))$ has a limit $T$, one has that $T$ contains the
tangent space at $x_0$ of $V_\a$, $T_{x_0}V_{\a} \subset T$.

\nn

Since the  Whitney stratification  $(V_\a)_{\a \in A}$ has finitely many strata containing
$0$ in its closure, essentially by \cite{Mi2} one has that each sufficiently small sphere
$S_{\varepsilon}$ intersects each stratum $V_\a$ transversally and the homeomorphism
type of the intersection
$K_{\varepsilon} = V \cap S_{\varepsilon}$ does not depend on $\e$ (c.f. \cite {LT,BV}). 
By Thom's transversality, if the stratification
satisfies the $a_f$-condition, then given such  $\e$ one has that for every
small disc $\D_{\eta}$ in $\R^k$ centred at $0$ and $t \in \D_{\eta}$,
  the level surface $f^{-1}(t)$ 
 also intersects transversally the sphere  $S_{\varepsilon}$. 
Let us set
$N(\e,\eta) = f^{-1}(\D_{\eta}) \cap \B_{\e}$, 
where $\D_{\eta}$ is a 
small disc in $\R^k$ centred at $0$ and $\B_{\e} $ is the ball bounded by
$S_{\varepsilon}$. Let 
 $T_{\varepsilon, \eta} = N(\e,\eta) \cap S_{\varepsilon} $, so it is
an algebraic neighbourhood of the link $K_{\varepsilon}$ in the sense of \cite {Du}.
 We refer to $N(\e,\eta)$ as a {\it Milnor tube} for $f$.

The following result is an extension of Milnor's fibration Theorem
 \cite[11.2]{Mi2}. This result is certainly known to various people and the first part of it is implicit in \cite {Le1},  but we did not find it  explicitly stated in the literature so we include it here with a short proof (which is essentially Milnor's original proof of his fibration theorem in \cite{Mi1}). 

\begin{theorem}{ Let the analytic map $f: X \to \R^k$ have $0 \in \R^k$ as  
an isolated critical value 
and assume there exists a  Whitney stratification $(V_\a)_{\a \in A}$ which 
satisfies the Thom $a_f$-condition with respect to $f$. Then 
for every $\e > 0$ sufficiently small and $\eta > 0$ 
sufficiently small with respect to $\e$, the map 
$$f: N(\e,\eta) \setminus V \to \D_{\eta} \setminus \{0\}\;,$$ is the projection of 
a locally trivial fibre bundle. This induces a fibre bundle projection:
$$\phi : S_{\varepsilon} \setminus Int(T_{\varepsilon, \eta})  \longrightarrow
S_{\eta}^{k-1} = \partial \D_{\eta}\;,$$ 
which restricted to the boundary 
$\partial T_{\varepsilon, \eta} = f^{-1}(S_{\eta}^{k-1}) \cap 
S_{\varepsilon} $ is the map ${f}$.}
\end{theorem}

\n
{\bf Proof.} Since  the stratification is Whitney, one may consider that there are only finitely many strata
having $0$ in its closure, and choose a small enough sphere $\S_\e$ around $0$ that intersects all strata transversally (and so does every other smaller sphere around $0$). 
 Since the Whitney stratification satisfies the $a_f$-condition, we can choose $\eta >0 $ 
sufficiently small so that for each
$t \in \D_{\eta} \setminus \{0\}$ the fibre $f^{-1}(t)$ intersects 
$S_{\varepsilon} $ transversally. Setting  
$N(\e,\eta) = f^{-1}(\D_{\eta}) \cap \B_{\e}$ as above, 
the first Thom-Mather isotopy theorem implies that the map:
$${f}:  N(\e,\eta) \to S_{\eta}^{k-1}\,, $$
is a locally trivial fibre bundle, proving the first claim in the theorem. This is also a direct application of Ehresman's fibration Theorem. 
 
Now let $f_1,...,f_k$ be the components of $f$ and define a function
$$r(x) = f_1^2(x) + ... + f_k^2(x)\;.$$
Let $\nabla r$ be its gradient. The level surfaces $r^{-1}(s)$, $s \in \R^+$, 
 are the tubes $f^{-1}(|t|)$ for $|t|^2 = s$ and the vector field 
 $\nabla r$ is transversal to these tubes (away from $V$) and tangent to $X$. 
Consider also the
vector field $\nabla(\iota)(x) = 2x$, which is the restriction to $X \setminus V$ of the 
gradient of the function
$x \buildrel{\iota} \over {\mapsto} x^2$. Both maps $r$ and $\iota$ 
are polynomial and $\ge 0$. Hence Corollary 3.4 of \cite{Mi2} implies that
on $\B_{\e} \setminus V$ the vector fields $\nabla r$ and $\nabla(\iota)$
are either linearly independent or one is a positive multiple of the other.
Following \cite{Mi1}, define a vector field on  $(X \cap \B_{\e}) \setminus V$ by:
$$v(x) = \|x\| \cdot \nabla r + \|\nabla r\| \cdot x\;,$$
which bisects the angle between $\nabla r(x)$ and $x$. This is a smooth 
vector field on $\B_{\e} \setminus V$ and satisfies the inequalities:
$$v \cdot x > 0, \hskip0,5cm v \cdot \nabla r > 0, \quad \ {(*)}$$
by the Schwartz inequality $\|v_1\|\|v_2\| \ge v_1 \cdot v_2$ (with equality 
holding only when the two vectors are colinear). 

Now, given any point $x_o \in (X \cap \B_{\e}) \setminus V$, let $\gamma(t)$ be the 
solution of the differential equation $dx/dt = v(x)$ in
 $ (X \cap \B_{\e}) \setminus V$ passing through $x_o$. As we move along the  path $\gamma(t)$, by $(*)$ 
the distance to $0$ increases strictly, and so does the function $r$, until 
$\gamma(t)$ intersects the boundary sphere. 

We may thus define a homeomorphism $h :  \partial N(\e,\eta) \to 
S_{\varepsilon} \setminus Int(T_{\varepsilon, \eta})$  as follows: given
$x \in \partial N(\e,\eta)$, take the path  $\gamma(t)$ passing through
$x$ and follow this path until it meets the sphere at a point $\gamma(t_1)$, 
then define $h(x) = \gamma(t_1)$. 
The composition $f \circ h^{-1}$ gives the fibration we want. \qed

In section 3 below we will see that the last part in the above proof can be refined when $f$
is holomorphic (Milnor's case) or, more generally, a product $h \overline g$ where $h$ and $g$ are both holomorphic mappings. This will lead to theorem 3.1.

\nn
{\bf 1.4 Examples.  a)}
Let  $f:(\C^2,0) \to (\C,0)$ and  $g:(\C^2,0) \to (\C,0)$ be two germs  of 
holomorphic functions  without 
common branches and such that $f \overline g : (\R^4,0) \to (\R^2,0)$ has an isolated 
critical value at $0$.
Let $f_1,\ldots,f_m$ be the irreducible factors of $f$ and let $g_1, \ldots, g_n$ 
be those of $g$. Let us equip a small open neighbourhood $U \subset \C^2$ of $0$ 
with the  Whitney stratification whose strata are :
$$ U \setminus (fg)^{-1}(0)\,;\, \; V_i = f_i^{-1}(0)\setminus \{0\}\,, \;
i=1,\ldots,m;$$ 
$$ V'_j= g_j^{-1}(0)\setminus \{0\}, j=1,\ldots,n ;\, \;\hbox{ and } \,  \{0\}\,.
$$ 
\n Then this stratification satisfies the $(a_f)$-condition with respect to 
$f \overline g $. In fact,
it suffices to check the condition on the 
strata $V_i, i=1,\ldots,m$ and $V'_j, j=1,\ldots,n$. Let us check this on 
$V_1$; the arguments are the same for the other strata.

Let us decompose $f$ as $f = f_1^p \,h$, $p \ge 1$, in such a way that $f_1$ 
is not a factor of $h$. Then the jacobian matrix of $f \overline g$ with respect to the 
coordinates $(z_1, \overline{z}_1, z_2, \overline{z}_2)$ in $\R^4$ is given by  
$$ D(f \overline{g})(z_1, \overline{z_1}, z_2, \overline{z_2})=
 \left( 
\begin{array}{cccc}
  \frac{\partial (\Re(f \overline g))}{\partial z_1} & 
\frac{\partial(\Re(f \overline g))}{\partial \overline z_1} & 
\frac{\partial(\Re(f \overline g))}{\partial z_2} & 
\frac{\partial(\Re(f \overline g))}{\partial \overline z_2}
  \\ 
  & & & \\
 \frac{\partial (\Im(f \overline g))}{\partial z_1} & 
\frac{\partial(\Im (f \overline g))}{\partial \overline z_1} & 
\frac{\partial(\Im (f \overline g))}{\partial z_2} & 
\frac{\partial(\Im (f \overline g))}{\partial \overline z_2} \\
\end{array}
\right)\;,
$$
where
$$\frac{\partial(\Re(f \overline{g}))}{\partial z_i} \, =\, 
\frac{1}{2} f_1^{p-1} \big( p\frac{\partial f_1}{\partial z_i}h \overline{g} 
+ f_1  \frac{\partial h}{\partial z_i}\overline{g} +\frac{\overline{f_1}^p}{f_1^{p-1}} 
\overline{h}  \frac{\partial g}{\partial z_i}  \big) \,:=\, \frac{1}{p} f_1^{p-1}m_{1i}\,,
$$
%\vskip.2cm
\[\frac{\partial(\Re(f \overline{g}))}{\partial \overline{z_i}} \, = \,
\frac{1}{p}\,\overline{ f_1^{p-1}m_{1i}}\;, \qquad  \qquad \qquad \qquad  \qquad \qquad \qquad  \qquad \qquad\]
%\vskip.2cm
$$\;\frac{\partial(\Im(f \overline{g}))}{\partial z_i} \, = \, \frac{-i}{2} f_1^{p-1} 
\big( p\frac{\partial f_1}{\partial z_i}h \overline{g} 
+ f_1  \frac{\partial h}{\partial z_i}\overline{g} -\frac{\overline{f_1}^p}{f_1^{p-1}} 
\overline{h}  
\frac{\partial g}{\partial z_i}  \big)
\, := \,\frac{1}{p} f_1^{p-1}m_{2i} \;,
$$
%\vskip.2cm
$$\frac{\partial(\Im(f \overline{g}))}{\partial \overline{z_i}}\, = \,
\frac{1}{p}\,\overline{ f_1^{p-1}m_{2i}} \;.\qquad  \qquad \qquad \qquad  \qquad \qquad \qquad  \qquad 
\qquad $$
\vskip.2cm
 For each $t \in \R^2 \setminus \{0\}$ and for each $x 
\in (f \overline g)^{-1}(t) $, the tangent space
$T_{x}(f \overline g)^{-1}(t)$ is defined by the equation: 

$$
D(f \overline{g})(x){}^t (v_1,\overline{v_1},v_2,\overline{v_2} )=0
$$
or, equivalently, by the equation 

$$M(x){}^t (v_1,\overline{v_1},v_2,\overline{v_2} )=0\,,$$
where: 

$$ M(x) = 
 \left( 
\begin{array}{cccc}
m_{11}(x) & \overline{m_{11}(x)} & 
 m_{12}(x) & \overline{m_{12}(x)} \\
  & & & \\
 m_{21}(x) & \overline{m_{21}(x)} & 
 m_{22}(x) & \overline{m_{22}(x)} \\
\end{array}
\right) \;.
$$

Let now  $\{x_n\}$ be a sequence of points in $R^4 \setminus \{0\}$ converging 
to a point $x \in f_1^{-1}(0)$. Then the limit $T$ of the tangent planes 
$T_{x_n}(f \overline g)^{-1}(t_n)$, where $t_n =(f \overline g)(x_n)$, has 
equation $M {}^t (v_1,\overline{v_1},v_2,\overline{v_2} )=0$, where 
$$M = \lim_{n\to\infty}  M(x_n)\,;$$

\noindent
but 
$$ \lim_{n \to \infty} m_{1i}(x_n) = \frac{1}{2} \, \frac{\partial f_1}{\partial z_i}(x) 
h(x) \overline{g(x)}\,,$$
and
$$ \lim_{n \to \infty} m_{2i}(x_n) = \frac{-i}{2} \, \frac{\partial f_1}{\partial z_i}(x) 
h(x) \overline{g(x)}\,.$$

As $h(x) \overline{g(x)} \ne 0$, then $T$ has equation:
$$N(x){}^t (v_1,\overline{v_1},v_2,\overline{v_2} )=0$$
where: 
$$ N(x) = \frac{1}{2} 
 \left( 
\begin{array}{cccc}
 \frac{\partial f_1}{\partial z_1}(x) & 
  \overline{\frac{\partial f_1}{\partial z_1}(x)}  &
 \frac{\partial f_1}{\partial z_2}(x)
 & 
  \overline{\frac{\partial f_1}{\partial z_2}(x)}  \\
   & & & \\
 -i\frac{\partial f_1}{\partial z_1}(x) & 
 -i  \overline{\frac{\partial f_1}{\partial z_1}(x)}  &
 -i \frac{\partial f_1}{\partial z_2}(x)
 & 
 -i \overline{\frac{\partial f_1}{\partial z_2}(x)} \\
\end{array}
\right) \;.
$$

Then $T$ equals the plane tangent 
 at $x$ to the curve $f_1^{-1}(0)$; so one has that these singularities satisfy the Thom $a_f$-condition for the above Whitney stratification.

\nn
{\bf Example b).}
Let now $(X,0)$ be a germ of a normal complex surface and let 
$f:(X,p) \to (\C,0)$ and  $g:(X,p) \to (\C,0)$ be two germs  of 
holomorphic functions  without 
common branches and such that $f \overline g : (X,p) \to (R^2,0)$ has 
$0$ as isolated critical value. Let us equip a small open neighbourhood $U \subset X$ of $p$ with the Whitney stratification defined as in example (a). Then this stratification  
 again satisfies the $(a_f)$-condition with respect to $f \overline g $. In fact, let 
$x \in X \setminus \{0\}$ and let $\phi : (W,0) \to (V,p)$ be a local chart, 
where $W$ and $V$ are open sets respectively  in $\C^2$ and $U \setminus \{0\} $. 
Then  the Whitney stratification of $U$  satisfies the $(a_f)$-condition with respect to $f \overline g $ if and only if its pull-back by $\phi$  satisfies the $(a_f)$-condition with respect to $\phi \circ f \overline g : (U,0)\to (\C,0)$. This last condition holds by similar 
arguments as in example (a).

\nn
{\bf Example c).}
 The previous examples can be generalized to higher dimensions. For instance consider  the maps in Example 1.1.c above. For simplicity restrict to $n =3$; one has:  
$$h(x,y,z) = \overline x \overline y \overline z (x^p + y^q + z^r)\,\,\,,\, p,q,r >1,$$
where the integers $p, q, r$ are $\ge 2$
and the unordered triple $(p,q,r)$ is not $(2,3,6)$, $(2,4,4)$ or $(3,3,3)$, so one has an isolated critical value at $0$.  Take the obvious 
stratification whose Whitney strata are $\C^3 \setminus h^{-1}(0)$; the coordinate planes minus the axes;
the complex surface $\{x^p + y^q + z^r = 0\}$ minus its intersection with the coordinate planes;  the 3 axes minus the origin; and  $\{0\}$. Then the computations of example a) above extend to this case and show that this stratification satisfies Thom's $a_f$-condition.

\nn
{\bf Example d).}
Let $f: (\C^2,0)\to (\C,0)$ and  $g: (\C^2,0)\to (\C,0)$ be two germs of holomorphic 
functions without common branches and such that 
$f \overline{g} : (\C^2,0) \to (\R^2,0)$ has $0$ as isolated critical value. 
Then for each holomorphic germ $h: (\C^n,0) \to (\C,0)$ with an isolated critical value
at $0 \in \C$  the real analytic germ 
$H: (\C^{n+2})\to (\R^2,0)$ defined by 
$$ H(z_1,z_2,z_3,\ldots,z_{n+2}) = f(z_1,z_2) \overline{g(z_1,z_2)}+ h(z_{3},\ldots,z_{n+2})
$$
obviously also has $0$ as an isolated critical value. Furthermore, by example (a) one has on 
$\C^2$ a Whitney stratification satisfying the Thom $a_f$ condition with respect to the function 
$f \overline{g}$. By \cite {Hi}  one also has a Whitney stratification of a neighbourhood of $0$ in $\C^n$
 satisfying the Thom condition for $h$ (see also \cite{LT}), and the spaces tangent to the fibres of $H$ split as direct 
 sum of the spaces tangent to the fibres of $f \overline{g}$ and $h$. Thus one also has the 
 Thom condition for $H$.

%\vskip1truecm
%\newpage

\section{A fibration theorem for functions $f \overline g$}

Let us consider now an equidimensional, reduced complex analytic isolated singularity $(X,0)$ of dimension 
$n$ in $\C^N$ and two germs of holomorphic functions $f,g: (X,0)\to(\C,0)$ 
such that the map $\, f \overline g : X \to \R^2\;,$
has an isolated critical value at $0$;
set $V = V_{f \overline g} := f \overline g^{-1}(0)$. 
Let us denote by $\mathcal{L}_X$ the link of $(X,0)$, {\it i.e.} 
the intersection of $X$ with a sphere $\S_\e$ around $0$ with radius $ \e >0$ 
sufficiently small, $\mathcal
 L_X = X \cap\S_\e $. We denote by
$L_{f \overline g}$ the link of $f \overline g$ in $\mathcal
 L_X$, {\it i.e.} the intersection of 
$V$ with 
$\mathcal
 L_X$ (for a possibly smaller $\e >0$). 

\begin{theorem} 
Let $(X,0)$ be a complex analytic isolated singularity germ and let
 $f,g$ be two germs of holomorphic functions $(X,0) \to (\C,0)$ such that 
$\, f \overline g : X \to \R^2\;$
has  $0$ as an isolated critical value. Then
$$ \phi = \frac{f \overline g}{|f \overline g|} : \mathcal{L}_X  
\setminus L_{f \overline g} \longrightarrow S^1\;,$$
is a $C^{\infty}$ fibre bundle. 
\end{theorem} 

\nn
 
 The proof of this theorem  takes the rest  of this section. The method we follow is exactly the same method used by Milnor in \cite[Chapter 4]{Mi2} to prove his fibration theorem for complex singularities: first show that $\phi$ has no critical points, so that all the fibres $\phi^{-1}(e^{i\t})$ are smooth manifolds of codimension 1 in the link
 $\mathcal{L}_X $;
  then construct a smooth vector field on the complement $\L _X  \setminus L_{f \overline g}$,
  whose integral lines are transversal to the fibres of $\phi$ and move at constant speed with respect
  to the argument of $\phi(z)$, so they carry each fibre diffeomorphically into the nearby fibres, 
  giving a local product structure.  We notice that this same method
is being used in \cite {BP} to prove a semitame condition for the meromomorphic function $f/g$ at the origin, {i.e.} a condition concerning the special values of $f/g$, and deduce from it a fibration theorem for these functions at 
infinity.

To begin  
we notice that away from the variety $ V_{f \overline g} = \{f \overline g = 0\}$ one has:
$f \overline g = \frac{f}{g} \,\cdot |g|^2$ and therefore:
$$\frac{f \overline g}{|f \overline g|} = \frac{f/g}{|f/g|}\;.$$
So to prove the theorem we actually show that this last map is the projection of a 
fibre bundle away from the link of $f \overline g$ in $\mathcal
 L_X$. 
The meromorphic function $f/g$ is well defined on 
$\mathcal{L}_X  \setminus L_{f \overline g} $ 
and takes values in ${\P}^1 \setminus \{0,\infty\}$.

Let us consider first the case when the germ $(X,0)$ is $(\C^n,0)$, so that 
$\mathcal  L_X = \S^{2n-1}$. 
Following \cite {Mi2}, define the gradient of the  meromorphic germ $h = f/g$ by : 
$$ \hbox{grad}\, h =(\overline{\frac{\partial h}{\partial z_1}}, \ldots, 
\overline{\frac{\partial h}{\partial z_n}})\, .$$

\begin{lemma} The critical points of the map 

 $$
\phi \, =\,\frac{f /g}{ \vert f/g  \vert} : {\S}^{2n-1}_{\epsilon}\setminus 
(L_f \cup L_g) \to {\S}^1\;,$$

 \noindent
are the points $z=(z_1,\ldots,z_n)$ where the vector 
$\big(i \,\hbox{grad} ( \log (\frac{f}{g}))\big )$ is a real multiple of $z$. 
\end{lemma}

\n{\bf Proof.} 
The proof of 2.2 is exactly as that of Lemma 4.1 in \cite {Mi2}, 
 setting $\frac{f /g}{ \vert f/g  \vert} = e^{i\theta(z)}$; we give it here for completeness and because some parts of the proof are also used later in the text. 
  Notice  that $\phi$ assigns to each $z \in \C^n$  the argument  
 of the complex number  $\frac{f} {g}(z)$:
 $$\phi(z) =   \frac{f /g}{ \vert f/g  \vert} (z) := e^{i\t(z)} \,,$$
 and therefore the argument of this map satisfies:
 $$\t(z) = Re \,(-i \,log \, f /g(z))\,.$$
An easy computation then shows that given any curve $z = p(t)$ in 
 $\C^n \setminus V_{f \overline g}$, the chain rule implies:
 \[d\t(p(t))/dt \,=\,  Re\,\langle \frac{dp}{dt}(t) , i \, grad\,log \,\frac{f}{g}(z) \rangle \,, \tag{2.2.1}\]
where $\langle \cdot, \cdot \rangle$ denotes the usual hermitian product in
$\C^n$.   Hence given a
   vector $v(z)$ in $\C^n $ 
 based at $z$, the directional derivative of $\t(z)$ in the direction of $v(z)$ is:  
 \[ Re\,\langle v(z), i \, grad\,log \,\frac{f}{g}(z) \rangle \,.  \]
  Since the real part of the hermitian product is the usual inner product in $\R^{2n}$, it follows that if
 $v(z)$ is tangent to the sphere ${\S}^{2n-1}_{\epsilon}$ then the corresponding directional derivative vanishes
 whenever  $\big(i \,\hbox{grad} ( \log (\frac{f}{g}))\big )$ is orthogonal to the sphere, {\it i.e.}  when it is a real multiple of $z$; conversely, if this inner  product vanishes for all vectors tangent  to the sphere then $z$ is a critical point of $\phi$. \qed
 
 \v
 
The aim now is to use lemma 2.2 to prove that $\phi$ does not have critical point at all.
This will be a consequence of the following lemma {(c. f. Lemma 4.4 in \cite {Mi2})}.

\begin{lemma} 
Let $p[0,\epsilon[ \to {\C}^n$ be a real analytic path with $p(0) =0$,
such that for each $t >0$ the number $f(p(t))g(p(t))$ is non-zero and 
the vector $\hbox{grad}(\log (f/g)(p(t)))$ is a complex multiple $\lambda(t) p(t)$. 
If $ \hbox{grad}\,(f /g) $ is non-zero away from the hypersurface $(fg)^{-1}(0)$ in a neighbourhood 
of the origin in $\C^n$,  then the 
argument of the complex number $\lambda(t)$ tends to $0$ as $t \to 0$.
\end{lemma} 

\nn
{\bf Proof.} 
One considers the Taylor expansions
$$p(t) = a t^{\alpha} + a_1 t^{\alpha +1} + a_2 t^{\alpha + 2} + \ldots\,;$$

$$f(p(t)) = b t^{\beta} + b_1 t^{\beta +1} + b_2 t^{\beta + 2} + \ldots\,;$$

$$g(p(t)) = b' t^{\beta'} + b'_1 t^{\beta' +1} + b'_2 t^{\beta' + 2} + \ldots\,;$$

$$\hbox{grad}(f(p(t))) = c t^{\gamma} + c_1 t^{\gamma +1} + c_2 t^{\gamma + 2} + \ldots\,;$$

$$\hbox{grad}(g(p(t))) = c' t^{\gamma'} + c'_1 t^{\gamma' +1} + c'_2 t^{\gamma' + 2} 
+ \ldots\,,$$
where $a, c, c' \in {\C}^n$ and $b,b' \in \C$ are non-zero. Notice that $\gamma 
= \beta - \alpha$ and   $\gamma' = \beta' - \alpha'$.

\n
Arguing as in p. 33 of \cite {Mi2} one gets:
$$\hbox{grad}(log((f/g)(p(t))) \,= \,\frac{1}{\overline {f(p(t))}} \ 
\hbox{grad}(f(p(t))) - \frac{1}{\overline {g(p(t))}} \ \hbox{grad}(g(p(t)))  \;;$$
\n
since  $\hbox{grad}(\log (f/g)(p(t))) = \lambda(t) p(t)$, we obtain: 
$$ {\overline {g(p(t))}}\, \hbox{grad}(f(p(t))) - 
 {\overline {f(p(t))}}\, \hbox{grad}(g(p(t))) = 
\lambda(t) p(t)  {\overline {f(p(t))}} {\overline {g(p(t))}}\,.$$ 
Let us consider  the non-zero term of lowest degree in the Taylor expansion of 
the left hand side of the equation above. 
It has the form $(\overline{b'} d - \overline{b} d')t^\delta $, where 
$\delta \geq \beta+\beta'-\alpha$. Substituting the Taylor expansions in the last
equality we get:
$$(\overline{b'}t^{\beta'}+\ldots)(ct^{\gamma}+\ldots) - 
(\overline{b}t^{\beta}+\ldots)(c't^{\gamma'}+\ldots)
= \lambda(t) (a t^{\alpha}+\ldots)(\overline{b}t^{\beta} + 
\ldots)( \overline{b'}t^{\beta'}+\ldots) ,$$
\noindent
that is,
\[(\overline{b'} d - \overline{b} d')t^\delta +\ldots = \lambda(t) 
(a \overline{b}\overline{b'}t^{\alpha + \beta+\beta'} + \ldots )\;. 
\]
Then, if $\lambda(t) = \lambda_0 t^{-\nu} + \ldots$ is the Laurent expansion 
of  $\lambda(t)$ (where 
$\nu = \delta - (\alpha + \beta+\beta')$), one has:
$$ \overline{b'}d-\overline{b}d' =  \lambda_0  a \overline{b}\overline{b'}\;,$$
which implies $\lambda_0 \ne 0$.  We now observe that the  formula of \cite[p. 33] {Mi2} for  
$h=f \overline g$ becomes:
\[\frac{d}{dt}\Big(\frac{f}{g}\big(p(t)\big)\Big) 
= \frac{1}{g(p(t))^2} \ \langle dp/dt, \ \overline{g} \ \hbox{grad}\,f 
- \overline{f} \ \hbox{grad}\, {g} \rangle\;.\]
Thus one has:
\[
g(p(t)) \ \langle\frac{d}{dt}, \hbox{grad}(f)\rangle   - f(p(t)) \ \langle\frac{d}{dt}, 
\hbox{grad}(g)\rangle ¥= \quad \qquad \qquad \qquad \qquad \qquad \quad \]
$ \qquad \qquad \qquad \qquad  \quad\qquad \qquad \qquad= \,g(p(t)) \frac{d}{dt}(f(p(t))) -f(p(t)) \frac{d}{dt}(g(p(t)))\,.$

\nn
Substituting the Taylor expansions in this equality we obtain:
$$\beta - \beta'= \vert a \vert^2 \alpha \overline{\lambda_0}\,.$$
\n
Therefore $\lambda_0 \in {\R}$. 
\qed

\nn

Before we continue with the proof of Theorem 2.1 we make sure that our function $f/g$
satisfies the hypothesis of Lemma 2.3, {\it i.e.} that its gradient 
is non-zero away from the hypersurface $(fg)^{-1}(0)$ in a neighbourhood 
of the origin in $\C^n$. This is granted by the following lemma.

\begin{lemma} If $f \overline{g} :({\C}^n,0) \to ({\C},0)$ has 
an isolated critical value at $0$, then $\hbox{grad} (f/g) $ is non-zero away from the set
$(fg)^{-1}(0)$ in a neighbourhood of $0$.
\end{lemma} 
\n
{\bf Proof.}
It suffices to prove that for each point $(z_1,\ldots,z_n) \in \C^n \setminus 
\{fg = 0\}$ with  
$\hbox{grad} (f/g)(z_1,\ldots,z_n) =0$, 
 the rank of the jacobian matrix of $f \overline{g}$ at $(z_1,\ldots,z_n)$ is not maximal. 
In other words, we must prove that if for each variable $z_i$, $i = 1,\ldots,n$, one has
(call it {\it equation} $(i)$):
$$  g \frac{\partial f}{\partial z_i} =  
f \frac{\partial g}{\partial z_i}  \;$$
\n
then all $2 \times 2$ minors of the jacobian matrix have zero determinant.
This is equivalent  to saying that
for each pair of variables $(z_i,z_j)$ one has the following equalities (c. f. \cite{Pi2}):
\[  fg \big( {\frac{\partial f}{\partial z_i}} \frac{\partial g}{\partial
z_j}  -  \frac{\partial f}{\partial z_j} \frac{\partial g}{\partial z_i}  \big) =  0 
 \qquad \quad  \qquad \qquad {(1.i.j)}\]

$$ \qquad \mid g \frac{\partial f}{\partial z_i}\mid   \,= \,\mid f  
\frac{\partial g}{\partial z_i} \vert  \quad   \qquad \;\qquad  \qquad \qquad     {(2.i.j)}  $$

$$ \qquad \mid g \frac{\partial f}{\partial z_j}\mid \, =  \,\mid f  
\frac{\partial g}{\partial z_j} \mid   \qquad   \quad \qquad \qquad \qquad    {(3.i.j)} $$

$$\mid g \mid^2  \frac{\partial f}{\partial z_i} 
\overline{\frac{\partial f}{\partial z_j}} \, = \,
\mid f \mid^2  \frac{\partial g}{\partial z_i} \overline{\frac{\partial g}{\partial z_j}} 
\ \   \qquad   \qquad \qquad  {(4.i.j)} $$

\n
Equations (2.i.j) and (3.i.j) are obvious from equation $(i)$ and the similar 
equation $(j)$: $  g \frac{\partial f}{\partial z_j} =  
f \frac{\partial g}{\partial z_j}  \;$. 
Now, multiplying the left hand side of equation $(i)$ by the right hand side of equation $(j)$ and vice-versa we obtain equation $(1.i.j)$. 
 Similarly, equation $(4.i.j)$ is obtained by conjugating equation $(j)$ and 
   multiplying left hand sides and right hand sides of the two equations.  \qed

\nn

From now on the proof of 2.1 is exactly like Milnor's proof.  By 2.2, 
to prove that $\phi$ has no critical points it is enough to show:

\begin{lemma} 
Given $f, g$   and  $V_{f \overline g}$ as above,
 there is a number $\e_o> 0$ such that for all $z \in \C^n \setminus V_{f \overline g}$
with $\Vert z \Vert \le \e_o$, the two vectors $z$  and $i \,grad\,log\, \frac{f}{g}(z)$ are either linearly
independent over $\C$ or else the argument of the complex number
$$\lambda \,=\, \big(i \,grad\,log\,\frac{f}{g}((z)) \big) / z$$
 has absolute value more than, say, $\pi/4$. 
 
 \end{lemma}
 
 The proof of 2.5 is exactly like that of Lemma 4.3 in Milnor's book, 
 using his Curve Selection Lemma and replacing \cite[4.4]{Mi2} by  lemmas 2.3 and 2.4 above. We leave the details to the reader. 
 
 \v
 As a consequence of 2.2 and 2.5 one gets that the map $\phi$ has no critical points at all.  
 It follows that all fibres of
 $\phi$ are smooth submanifolds of $L_{f \overline g}$ of real dimension $2n-2$. 
In order to show that $\phi$ is actually the projection map of a $C^\infty$ fibre bundle we must prove that one has a local product structure around each fibre. As in \cite{Mi2}, this is achieved by making a sharper use of 2.5 to construct a vector field $v$ on 
${\S}^{2n-1}_{\epsilon}\setminus (L_f \cup L_g) $ 
so that  at each point $z \in {\S}^{2n-1}_{\epsilon}\setminus (L_f \cup L_g) $ 
one has that the hermitian product  
\[\langle v(z), i\,grad \,log\,\frac{f}{g}(z) \rangle \,,\]
is non-zero and has argument less than $\pi/4$ in absolute value; in particular its real part is never-zero.
Thus we can normalize it and define:
\[w(z) \,=\, v(z) \big / Re\,\langle v(z), i\,grad \, log \,\frac{f}{g}(z) \rangle \;.
\] 
 This is a $C^\infty$ vector field on ${\S}^{2n-1}_{\epsilon}\setminus (L_f \cup L_g) $ satisfying:
 
 \nn
 i) the real part of the hermitian product $\langle w(z), i\,grad \, log \,\frac{f}{g}(z) \rangle$ is identically  equal to 1;  recall that by equation (2.2.1) this is the directional derivative of
 the argument of $\phi$ in the direction of $w(z)$.

 \nn
 ii) the absolute value of the corresponding imaginary part is less than 1:
 \[\vert Re\, \langle w(z), \,grad \, log \,\frac{f}{g}(z) \rangle \vert \,< \,1\,. \tag{2.6}\]
  
 \v
 Let us consider now the integral curves of this vector field, {\it i.e.} the solutions $p(t)$  
 of the differential equation
 $dz/dt = w(t)$. Set $e^{i\t(z)} = \phi(z)$ as before.  Since the directional derivative of $\t(z)$ in the direction $w(t)$ is identically equal to 1 we have:
 $$\t(p(t)) \,=\, t + constant\,.$$
 Therefore the path $p(t)$ projects to a path which winds around the unit circle in the positive direction with unit velocity. In other words, these paths are transversal to the fibres of $\phi$ and for each $t$ they
 carry a point
 $z \in \phi^{-1}(e^{it_o})$ into a point in $ \phi^{-1}(e^{it_o+t})$. If there is a real number $t_o > 0$ so that all
 these paths are defined  for, at least, a time $ t_o$ then,   being solutions of the above differential equation, they will carry each fibre of $\phi$ diffeomorphically into all the nearby  fibres, proving that one has a local product structure and $\phi$ is the projection of a locally trivial fibre bundle. In fact, just as in 
 Lemma 4.7   of  \cite{Mi2},
 one has that all these paths extend for all $t \in \R$.
To show this  we only have to   insure that no path $p(t)$ can converge to $(L_f \cup L_g)  $ as $t$ tends to some finite time. That is, we must show that $f/g(p(t))$  can not tend to 0 in  finite time, or equivalently that the real part of $\,log\,f/g(p(t))\,$ can not tend to $-\infty$ as $t$ tends to some finite time 
$t_o$, but this is an immediate consequence of 2.6.
  Thus we arrive to Theorem 2.1  when 
$X= \C^n,0$. 

The general case of $X$ being a complex analytic variety with an isolated singularity
 follows from the previous considerations since  all the above arguments are local. \qed

\section{The two fibrations are equivalent}

Let $(X,0)$ be again an  equidimensional, reduced
 complex analytic isolated singularity germ of dimension $n$ in $\C^N$,
 let
 $f,g$ be two germs of holomorphic functions $(X,0) \to (\C,0)$ such that 
$$ f \overline g : X \to \R^2\,,$$
has  $0$ as an isolated critical value. We believe one always 
 has in this situation a Whitney stratification of $X$ that satisfies 
the Thom $a_f$-condition with respect to $f\overline{g}$, but we do not know how to prove this statement (example 1.1.b above proves it for $n = 2$).  
Anyhow we show below that 
 one always has a fibre bundle, 
\[f \overline g : N(\e,\eta)  \to \partial \D_{\eta}   \,, \]
as in Theorem 1.1, for every sufficiently  small circle
$\partial \D_{\eta}   \subset \R^2$ 
 of radius $\eta$ around the origin in $\R^2$, where $N(\e,\eta)$ 
is the Milnor tube 
$N(\e,\eta) = (f \overline g)^{-1}(\partial \D_{\eta}) \cap \B_{\e}$. We recall
one has the  fibre bundle given by Theorem 2.1:
$$ \phi = \frac{f \overline g}{|f \overline g|} : \mathcal{L}_X  \setminus L_{f \overline g} \longrightarrow S^1\;.$$
Let $T(\e,\eta)$ be the open (regular) neighbourhood of the link $L_{f \overline g}$ in $\L_X$
bounded by  $N(\e,\eta) \cap \L_X$.  One has:

\begin{theorem}  These two fibrations are topologically equivalent. More precisely,  there exists a $C^\infty$ 
diffeomorphism $h: N(\e,\eta) \setminus V_{f \overline g} \to \mathcal{L}_X  \setminus  T(\e,\eta)$ such that $f|_{N(\e,\eta)} \,= \, \phi \circ h$; and $h$ is the identity on $N(\e,\eta) \cap \L_X$. 
\end{theorem} 

In fact the proof of this result is similar to that of Theorem 1.1 above but in reverse sense: in 1.1 we first used Thom's condition and Ehresman's lemma to deduce the fibre bundle structure on the Milnor tube, and then we constructed an appropriate vector field to inflate the tube and get a fibration on the  link $\mathcal{L}_X$ of $X$ minus the  tube $T(\e,\eta)$.
 Now we have Theorem 2.1 that gives a fibre bundle decomposition on 
$\mathcal{L}_X  \setminus L_{f \overline g} $ with projection map 
$\phi = \frac{f \overline g}{|f \overline g|} \, =\, \frac{f/g}{\vert f/g \vert}$; then we construct again an appropriate vector field that 
carries this fibre bundle into the Milnor tube 
$(f/g)^{-1}(\partial \D_{\eta}) \cap \B_\e$,
 and whose paths move along points in $X$ where the argument of the complex number  $f/g (z)$ does not change. 
\v

The key for proving this theorem is the following lemma (c.f. 5.9 in \cite {Mi2}):

\begin{lemma}  Let  $\D_\e$ be the intersection with $X$ of  a sufficiently small ball in $\C^N$ around $0$. Then there exists a $C^\infty$ vector field $v$ on $ \D_\e \setminus V_{f \overline g}$ so that for all 
$z \in \D_\e \setminus V_{f \overline g}$ one has:

\nn {\bf i)}  the hermitian product 
$ \langle v(z) , \, grad\, log\, \frac{f}{g}(z) \rangle $
is real and positive;  and 

\nn {\bf ii)}   the hermitian product  $\langle v(z), z \rangle $ has positive real part.
\end{lemma}

\n
{\bf Proof.}   It suffices to construct this vector field locally, in the neighbourhood of some given point $z_o \in \D_\e \setminus V_{f \overline g}$, for we can afterwards glue all these vector fields by a partition of unity and get a vector field defined globally and satisfying the conditions of the lemma. We work on a coordinate chart for $X$ at $z_o$, so we think of it as being identified with an open ball $U$ in $\C^n$.
We recall that the real part of the hermitian product is the usual inner product in $\R^{2n}$. 

For simplicity set $\nabla(z_o) = grad\, log\, \frac{f}{g}(z_o)$,
 and let
$E$ denote the real orthogonal complement of $ i \nabla(z_o)$ in $\C^n$; every  $w$ in $E$ satisfies   
that $ \langle v(z_o) , \nabla(z_o)\rangle $ is real; the vector  $  \nabla(z_o)$  defines an oriented real line $\l$ in $E$ and  every vector in the half-space $E^+$ of vectors in 
$E$  whose projection to $\l$ is positive satisfies condition i) in the lemma. We need to show that there
are vectors in $E^+$ satisfying also condition ii) in 3.2. If $z_o$ and $\nabla(z_o)$ are linearly dependent over $\C$, so that 
 $\lambda z_o = \nabla(z_o)$  for some $\lambda$ in $\C$, then 2.5 above implies:
$$\vert argument \, \lambda \vert < \pi/4\,,$$
and we can just take $v = z_o$.
If the vectors $z_o$ and $\nabla(z_o)$ are linearly independent over $\C$, then
  $  \nabla(z_o)$  and $z_o$ are linearly independent over $\R$ and every vector in 
  $E^+$  for which the real inner product with $z_o$ is positive works for us. \qed

\nn

We now return to the proof of 3.1. Consider the solutions of the differential equation
$$dp(t) \,=\, v(p(t))$$
on $\D_\e  \setminus V_{f \overline g}$. By the first statement in 3.2 we have that the argument of 
$f/g(p(t))$ is constant, and that $\vert f(p(t))$ is strictly monotone as a function of $t$.  The second statement in 3.2 insures that $\Vert p(t)\Vert$ is a strictly monotone function of $t$. Thus starting at any
interior point $z$ of $\D_\e  \setminus V_{f \overline g}$ we travel away from the origin, in a direction of increasing
$\vert f\vert$ and constant argument, until we reach a point 
$z' \in \mathcal{L}_X  \setminus L_{f \overline g} $ with 
$$\frac{f/g(z)}{\vert f/g(z)\vert} \,=\, \frac{f/g(z')}{\vert f/g(z')\vert}\,.$$
Let $\eta >0 $ be small enough,  as in Theorem 3.1, and look at the Milnor  tube 
$$N(\e,\eta) \,=\, \{z \in \D_\e  \setminus V_{f \overline g} \, \big \vert \, f/g(z) \,= \eta\,. \}$$
Then the  correspondence $z \buildrel {\hat h} \over \mapsto z'$ carries this tube into the complement of a tubular neighbourhood of $ L_{f \overline g} $ in 
$\mathcal{L}_X$. One obviously has:
$$ {f/g} \,=\, \phi \,\circ\, \hat h \,, $$
for all $z$ in the tube $N(\e,\eta)$. This identifies the complement of a 
tubular neighbourhood of
of $L_{f \overline g}$ in $\mathcal{L}_X$ with the Milnor tube for
$f/g$. To complete the proof of 3.1 we now carry this tube diffeomorphically
 into a Milnor tube of $f \overline g$. 
 For this we notice that for all $z$ in $\D_\e  \setminus V_{f \overline g}$ the map
$ z \mapsto \frac{z}{\vert g(z) \vert^2}$ is a well defined diffeomorphism that we denote by $Q$. Then
one has 
$f \overline g \,=\, {f/g} \, \circ \, Q\,, $ and therefore
$$ f \overline g \,=\,    \phi \,\circ\, \hat h \,    \, \circ \, Q\,, $$
 proving 3.1 with $h = \hat  h \circ  Q$. \qed

 \nn
 \begin{corollary} The fibres $F_\t = \phi^{-1}(e^{i \t})$ are diffeomorphic to the complex $(n-1)$-manifold consisting of the intersection $f/g^{-1}(t) \cap \buildrel {\circ} \over {\D}_\e$, where $t$ is a regular value of the rational map $f/g$ and $\buildrel {\circ} \over {\D}_\e$ is the interior of the disc of radius 
 $\e >0$ around the origin.
   Hence each fibre $F_\t$ has the homotopy type of a CW-complex of real dimension $(n-1)$.
\end{corollary}

\n
{\bf Proof.}  The proof of 3.1 shows that $f/g^{-1}(t) \cap \buildrel {\circ} \over {\D}_\e$ is diffeomorphic to
the part of $F_\t = \phi^{-1}(e^{i \t})$ contained in the complement of a regular neighbourhood of the link
$L_{f \overline g}$, so the theorem of Andreotti-Frankel in \cite {AF} implies that this part of
 $F_\t$ has the homotopy type of a CW-complex of real dimension $(n-1)$.
We claim that this part of $F_\t $ is actually diffeomorphic to the whole
 fibre $F_\t$.  The proof of this claim is similar to that of \cite[5.7]{Mi2}:  it is enough to show that the real valued map $a_\t: F_\t  \to \R$ given by $z \mapsto \vert f/g(z) \vert $ has no critical values near the binding $ L_{f \overline g} $, or more precisely that there exists a constant $\delta_\t > 0$ such that all the critical points of $a_\t$ lie within the compact subset 
 $\vert f \overline g(z) \vert  \ge  \delta_\t $ of  $F_\t$. Suppose this does not happen; then one must have a sequence of critical points of $a_\t$ which converges to some limit point 
 $z_o \in L_{f \overline g} \subset
 \partial \D_\e$. By the Curve Selection Lemma of \cite {Mi2} one  has a curve $\gamma$ of
 points in $F_\t$ converging to $z_o$, consisting of critical points of $\a_\t$. Clearly $a_\t(z)$ is constant
 along $\gamma$ and can not converge to $0 = \vert f \overline g (z_o) \vert$. \qed

\nn
 \begin{corollary} If the tangent bundle $T(X \setminus\{0\})$ of the ambient space is trivial, 
 either as a real or a complex bundle, then the tangent bundle of each fibre $F_\t = \phi^{-1}(e^{i \t})$ 
 is  trivial, $C^\infty$ isomorphic to either $F_\t \times \R^{2n-2}$ or to $F_\t \times \C^{n-1}$,  respectively.
\end{corollary}

\n
{\bf Proof.}  The previous result
implies that the homology groups  of $F_\t $ vanish in dimensions more than $n-1$, so $F_\t$ can not have any compact component. On the other hand, the normal bundle $\nu \hat F_\t$
of $\hat F_\t  =  (f/g)^{-1}(e^{i\t}) \cap \D_\e$ in $X$ is trivial over $\C$, because $t$ is a regular value. Since one has a splitting (as $C^\infty$ bundles)
$$T(X \setminus\{0\}) \vert_{F_\t} \,\cong \, TF_\t  \oplus \nu \hat F_\t \,,$$
it follows that  the tangent bundle  $TF_\t $ plus a trivial bundle is trivial (either over $\R$ or over $\C$) by hypothesis, {\it i.e.} $TF_\t $  is stably trivial. 
Hence it is trivial by \cite {KM}, because $F_\t$ has no non-compact components.

\section{Fibrations of plumbing multilinks }

From now on the standard circle $\mathbf S^1 = \{ z \in \mathbf C \ / \ |z| =1\}$ is oriented as the boundary of the complex disk $\mathbf D^2 = \{ z \in \mathbf C \ / \ |z| \leq 1\}$. A {\it circle} means a $1$-dimensional manifold diffeomorphic to 
$\mathbf S^1 $. 
 
 Let $M$ be a compact connected oriented $3$-manifold. An (unoriented) {\it knot} in $M$ is a circle embedded in $M$. An (unoriented) {\it link} $L$ in $M$ is a finite disjoint union of circles embedded in $M$. We denote $L = K_1 \cup \ldots \cup K_l$, where the $K_i$'s are the connected components of $L$.
 We say that $L$ is an {\it oriented link} if an orientation is fixed on each of its components $K_i$. Since the manifold $M$ is oriented, the normal bundle of each 
$K_i$ in $M$ is trivial  and $K_i$ admits a regular neighbourhood $N(K_i)$ diffeomorphic to the solid torus ${\mathbf D^2} \times {\mathbf S^1}$ by an orientation-preserving diffeomorphism which sends the oriented circle $K_i$ to $\{0\} \times \mathbf S^1$. 
 
 A {\it multilink } is the data of  an oriented link $L = K_1 \cup \ldots \cup K_l$ together with a multiplicity $n_i \in {\mathbf Z}$ associated with each component $K_i$. We denote such a multilink by 
$$L = n_1 K_1 \cup \ldots \cup n_l K_l\,,$$ and we fix the convention that $n_i K_i = (-n_i) (-K_i)$, where $-K_i$ means $K_i$ with the opposit orientation. We say that two multilinks $L$ and $L'$ are {\it isotopic} if they are isotopic as oriented links in $M$ and if their multiplicities are the same through such an isotopy.

\nn
{\bf 4.1 Example.} Let $(X,p) $ be a normal complex surface singularity and let $(X,p) \subset ({\mathbf C}^N,0)$ be an embedding. Let $\L_X$ be the link of $(X,p)$, {\it i.e.} the oriented compact connected $3$-manifold $X \cap \S_\e^{2N-1}$,
 intersection of 
$X$ and the sphere $\mathbf S^{2N-1}_{\epsilon} = \{ (z_1,\ldots,z_N) \in {\mathbf C}^N \ / \ \sum_{i=1}^N |z_i|^2 = \epsilon^2 \}$, $\epsilon >0$ sufficiently small (\cite{Mi2}). Let $f:(X,p) \longrightarrow
({\mathbf C},0)$ be a holomorphic germ, and let $L_f = \L_X \cap f^{-1}(0)$ be its link in $\L_X$, which is naturally oriented as the boundary of the complex curve $ f^{-1}(0) \cap {\mathbf D^{2N}_{\epsilon}}$, where ${\mathbf D}^{2N}_{\epsilon}=\{ (z_1,\ldots,z_N) \in {\mathbf C}^N \ / \ \sum_{i=1}^N |z_i|^2 \leq \epsilon^2 \}$.

Let $\pi : \tilde{X} \to X$ be a resolution of the germ $f$, {\it i.e.} a resolution of $X$ such that  $(f \circ \pi)^{-1}(0)$ is a normal crossings divisor. Let $S_1,\ldots,S_l$ be the branches of the strict transform $\overline{\pi^{-1}(f^{-1}(0) \setminus \{p\})}$ of $f$. We denote by $n_i$ the multiplicity of $f \circ \pi$ along the curve $S_i$. For each $i=1,\ldots,l$, let $K_i$ be the component  of $L_f$ associated with $S_i$, {\it i.e.} the knot $\pi(S_i \setminus \pi^{-1}(p)) \cap \L_X$. By definition  {\it the multilink associated with} $f$ is:
$$L_f  = \bigcup_{i=1}^l n_i K_i\,.$$

In particular, when the structural sheaf ${\cal O}_X$ is a unique factorization ring (e.g. $X={\mathbf C}^2$), let $f=\prod_{i=1}^l f_i^{n_i}$ be the  decomposition of $f$ into irreducible analytic factors, where $u(0) \not= 0$, $f_1,\ldots,f_l$ are irreducible, and $n_i \in \mathbf N^{\ast}$. Then  the multilink associated with $f$ is 
$$L_f  = \bigcup_{i=1}^l n_i L_{f_i}\,.$$

\nn
{\bf 4.2 Definition.}
Let $M$ be a compact connected oriented $3$-manifold. A multilink $L = n_1 K_1 \cup \ldots \cup n_l K_l$ in  $M$ is {\it fibred} if there exists a map $\Phi : M \setminus L \longrightarrow {\mathbf S^1}$ which satisfies the two following two conditions : 
\begin{enumerate}
\item The map $\Phi $ is a $C^{\infty}$ locally trivial fibration.
\item For each $i =1,\ldots,l$, there exist a regular neighbourhood $N(K_i)$ of $K_i$ in $M \setminus (L \setminus K_i)$, an orientation-preserving diffeomorphism  $\tau : {\mathbf S}^1 \times {\mathbf  D}^2
\to N(K_i)$ such that $\tau ({\mathbf S}^1 \times \{0\} )= K_i$ and an integer $k_i \in {\mathbf Z}$ such that for all $(t,z) \in \mathbf{S}^1
\times ({\mathbf  D}^2 \setminus \{0\})$ one has:
$$ (\Phi \circ \tau)(t,z)=
\bigg( \frac{z}{\vert z \vert} \bigg)^{n_i} t^{k_i} $$
\end{enumerate}
In this case we say that $\Phi$ is a {\it fibration of the multilink} $L$.
 
 \nn{\bf 4.3 Example.} 
If $(X,p)$ is a normal complex surface singularity and if $f:(X,p) \longrightarrow
({\mathbf C},0)$ is a holomorphic germ, then the Milnor fibration 
$$\Phi_f : \L_X \setminus L_f \to {\mathbf S}^1\,,$$
 defined by  $\Phi_f(x) = f(x) \ / \ |f(x)|$, is a fibration of the multilink $L_f$.

\nn {\bf 4.4 Remarks:}

\v \n {\bf  1)}  The integer $k_i$ is non unique ;  it depends on the choice of $\tau$, but its class modulo $n_i$ is well defined.

\v \n {\bf   2)} For each $t \in {\mathbf S}^1$,  $\Phi^{-1}(t) \cap (N(K_i) )$ has $gcd(n_i,k_i)$ connected components, each of them being diffeomorphic to a half-open annulus ${\mathbf S}^1 \times [0,1[$.

\v \n {\bf  3)} The $n_i$ are local date, whereas the $k_i$ are global. Indeed,  let $D$ be a meridian disk of $N(K_i)$ oriented as $\{1\} \times {\mathbf D}^2$ via $\tau$ and let us equip its boundary $m = \partial D$ with the induced orientation. Then $n_i$ is the degree of
the restriction of $\Phi$ to $m$. But $k_i$ depends on the isotopy class of the multilink $L$ in $M$, and in particular on all the multiplicities $n_i, i=1,\ldots,n_l$. See \cite{EN}  page 30,  where the inverse $l_i$ of $k_i$ modulo $n_i$ is
 computed explicitly when the link is an integral homology sphere. 

\v \n {\bf  4)}  If for each $i=1,\ldots,n$, $k_i=0$  and $n_i  \in \{-1,+1\}$, 
 one obtains respectively 
 $$ (\Phi \circ \tau)(t,z)=
\frac{z}{\vert z \vert}  \  \ \hbox{ when }  n_i = +1\,, \hbox{ or }$$
$$ (\Phi \circ \tau)(t,z)=
\frac{\bar{z}}{\vert z \vert}   \  \ \hbox{ when }  n_i = -1 \,. \ \  \ \ $$ 
In this particular case, which is studied in \cite {Pi2} when $M = {\mathbf S}^3$, $L$ is just an oriented link and $\Phi$ is a so-called open-book fibration of $M$ with binding $L$.

 \nn

Now let  $M$ be a {\it plumbing manifold}, or equivalently, a graph manifold in the sense of Waldhausen. In other words, $M$ is homeomorphic to the boundary of a $4$ manifold obtained by plumbing together a finite number of $\mathbf D^2$-bundles $N_1,\ldots,N_r$ over some oriented real surfaces $E_1,\ldots,E_r$ (see \cite{N}). 
If $T \subset M$ denotes the union of the plumbing tori, then $M \setminus T$ is a union of $3$-manifolds equiped with an $\mathbf S^1$-foliation, the $\mathbf S^1$-leaves being boundaries of the $\mathbf D^2$-fibres of the $N_i$'s.

\v
\noindent
{\bf Convention.} In the sequel we always assume that the surfaces $E_i$ are orientable and that an orientation on each $E_i$ is fixed. We also assume that each plumbing operation has positive $\epsilon$ in the sense of \cite{N}.

The homeomorphism class of a plumbing manifold $M$ is completely encoded in the  {\it plumbing graph} $\Gamma$ of any plumbing decomposition of $M $ defined by the classical way: the vertices $(1),\ldots, (r)$ are in bijection with the bundles $N_1, \ldots , N_r$, and the edges joining two vertices $(i)$ and $(j)$ are in bijection with the plumbing operations between $N_i$ and $N_j$. Each vertex $(i)$ is weighted by the genus of the surface $E_i$ (we usually omit to write this weight when the genus is $0$) 
 and by the Euler class $e_i$ of the bundle $N_i$. 
 
 We denote by  $M_{\Gamma}$ the {\it intersection matrix} associated with $\Gamma$, {\it i.e.} the $r \times r$ matrix $(a_{ij})_{1 \leq i,j\leq r}$ defined by  $ a_{ii} = e_i$  for all i and  for  $i \ne j$ the integer 
$a_{ij}$ equals the number of edges of $\Gamma$ joining the vertices $(i)$ and $(j)$.

\nn {\bf 4.5 Example.} (c.f. \cite{N}) The link $\L_X$ of a  normal complex surface $(X,p)$ is a plumbing manifold  and a plumbing graph of $\L_X$ can be obtained as the dual graph $\Gamma$ of any resolution $\pi: \tilde{X} \to X$ of $X$ whose exceptional divisor $E=\pi^{-1}(p)$  has normal crossings.  Let $E_1,\ldots,E_r$ be the irreducible components of $E$. The intersection matrix $M_{\Gamma}$ is then  the intersection matrix of the $E_i$'s in $\tilde{X}$, namely
$$M_{\Gamma} = \big( E_i.E_j\big)_{1 \leq i,j\leq r}\;.$$
The matrix $M_{\Gamma}$ is negative definite \cite{G}.

\nn{\bf 4.6 Definition.} Let $M$ be a plumbing manifold. A  link $L$ in  $M$ is a {\it plumbing link} if it is a finite union of $\mathbf{S}^1$-leaves  in a plumbing decomposition of  $M$, or equivalently,  a finite union of Seifert fibres in a Waldhausen
decomposition of  $M$.
 
 \v
Notice that any  plumbing link $L $ in $M$ has a natural orientation as the union of the boundaries of some fibres of the $D^2$-bundles $N_i$. In the sequel we consider
 always  this natural orientation. 

The  isotopy class of a plumbing link $L$ is encoded in a plumbing graph obtained by decorating a plumbing graph of $M$  with arrows corresponding to the components of $L$. Similarly, if $L = n_1 K_1 \cup \ldots \cup n_l K_l$ is a plumbing multilink, then its isotopy class is encoded in its plumbing graph, obtained from the previous decorated graph by weighting  the extremity of each arrow  by the corresponding multiliplicity $n_i$.

For example if $(X,p) $ is a  normal complex surface singularity and  $f:(X,p) \to
({\mathbf C},0)$ is a holomorphic germ, then $L_f$ is a plumbing multilink and a plumbing graph of $L_f$ can be obtained as  the dual graph of any resolution $\pi: \tilde{X} \to X$ of $f$ such that the total transform of $f$ by $\pi$ has normal crossings.

\v

\noindent
{\bf 4.7 Definition.}
Let $\Gamma$ be a plumbing graph of a plumbing link. A vertex of $\Gamma$ is a {\it rupture vertex} if it
either  carries genus $>0$ or if it has at least three incident edges, the arrows being considering as edges.
 
\v
The following theorem is a generalization of  a  result of D. Eisenbud and W. Neumann (\cite{EN}, Theorem 11.2, see also p. 136),
reformulated in terms of plumbing links. In \cite{EN}, this result is for  graph multilinks in  $\mathbf Z$-homology spheres, and it is formulated in terms of splicing diagrams. In $\cite {Pi1}$ this result is proved for a plumbing link in terms of Waldhausen graphs; here we extend this for multilinks. The proof follows that of  $\cite {Pi1}$, interpreting a plumbing decomposition of $M$  as a particular Waldhausen decomposition. 

\begin{theorem}
\label{fibration1}
 Let $L$ be a plumbing multilink with plumbing graph $\Gamma$ and intersection matrix $M_{\Gamma}$, and let $(1),\ldots,(r)$ be the vertices of $\Gamma$. Let $b(L)=(b_1,\ldots,b_r) \in {\mathbf{Z}^{r}}$, where $b_i$ is the sum of the multiplicities $n_i$ carried by  the arrows attached to the vertex $(i)$.
Let  $(m_1,\ldots,m_r) \in \mathbf{Q}^r$ be defined by 
$$ {}^{t}(m_1,\ldots,m_r) =-(M_{\Gamma})^{-1}\hspace{1mm} {}^{t} b(L)$$
where ${}^t(.)$ means the transposition.
Then $L$ is fibred if and only if the following two conditions hold:

1) for all $i \in \{1,\ldots,r\}$, $m_i$ is an integer; and

2) for each rupture vertex $(j)$ of $\Gamma$, the integer $m_j$ is  $\neq 0$.
\end{theorem}

\begin{proof} 
 Assume that $L$ is fibred. Let $ \Phi : M \setminus L \to {\mathbf S}^1$ be a fibration of $L$ and let $F = \Phi^{-1}(t), \ t \in {\mathbf S}^1$ be a fibre of $\Phi$. Notice that $F$ is naturally oriented by the flow. For each $i = 1, \ldots ,r$, let $l_i$ be a ${\mathbf S}^1$-leave of $N_i$. 
Then $m_i$ is the intersection number $F.l_i$ in $M$ (\cite{Pi1}). Moreover, when $\Gamma$ has at least one rupture vertex, then according to \cite{Wa}, $F$ is a horizontal surface in $M \setminus N(L)$, where $N(L)$ is a regular neighbourhood of $L$ in $M$. It means that $F$ is transversal up to isotopy to $L_j$ for each rupture vertex $(j)$ of $\Gamma$, {\it i.e.} $m_j \neq 0$.

Assume now that  $ N= (m_1,\ldots, m_r)$ verifies the conditions 1. and 2. Using the formulas of \cite{Pi1}, one constructs  a graph $G^{N}$ from $\Gamma$ and $N$. The  conditions 1. and 2. imply that $G^N$ is the Nielsen graph of the monodromy of a fibration of the multilink $L$ (\cite{Pi1}, lemma 4.7).  Therefore $L$ is fibred.  \qed
\end{proof}

 \nn {\bf 4.8 Remarks:}

\v \n {\bf 1)} An easy computation shows that each transformation $R_0$ to $R_{10}$ of (\cite{N})  which transforms a plumbing graph of $L$ into another one preserves  the conditions 1 and 2 in 4.1. Hence these conditions hold for some plumbing graph $\Gamma$ of $L$ iff  they  hold for every plumbing graph of $L$.
 
\v \n {\bf 2)} Given any plumbing multilink 
$L= n_1 K_1 \cup \ldots \cup n_l K_l$ one has that
the condition 1 holds always for the multilink $kL = (k n_1)  K_1 \cup \ldots \cup (k n_l) K_l $, where $k=\det M_{\Gamma}$. Thus, if condition 2 holds for the multilink $L$, then the multilink $k L$ is fibred.

\v \n {\bf  3)} Let  $L $ and $L '$ be two fibred multilinks respectively in $M$ and $M'$. We say that two fibrations $\Phi : M \setminus L \to {\mathbf S}^1$ and  $\Phi' : M' \setminus L' \to {\mathbf S}^1$ of $L$ and $L'$ respectively are {\it topologically equivalent} if there exist two orientation-preserving diffeomorphisms $H:(M,L)
\longrightarrow (M',L')$ and $\rho: {\mathbf S}^1 \longrightarrow
{\mathbf S}^1$
such that:
 \v {\bf i)} $\rho \circ \Phi = \Phi' \circ H{\mid_{(M \setminus
L)}}$;
\v {\bf ii)}  for each component $K$ of $L$, the multiplicity of $K$ in $L$ equals that of $H(K)$ in $L'$.
 
\n Note that this implies  in particular that the two multilinks are isotopic when $M=M'$. We now let $L$ be a fibred plumbing multilink in a plumbing manifold $M$. Then an immediate generalization of  the  arguments of \cite{Pi1} to multilinks shows that there exists a finite number of topological equivalence
classes  of fibrations of $L$. 
 
 \noindent
\nn {\bf 4.9 Examples:}
1) Let $M$ be the plumbing 3-manifold whose graph is $-D_4$ (see figure 1), and let $L_1,L_2,L_3$ and $L_4$ be the four multilinks whose graphs are represented on figure 1.  The numbers without  parenthesis are the Euler classes and the numbers with parenthesis are the multiplicities $m_i$ and $n_i$. Then $L_1$ is not fibred as the rupture vertex carries multiplicity $0$. The multilink $L_2 = L_3 -L_4$ is fibred since it verifies the conditions 1 and 2 of \ref{fibration1}, whereas $L_3$ and $L_4$ are not fibred since they do not verify the condition 1 of \ref{fibration1}. 

\v
2) Let $(X,0)$ be the germ of hypersurface in $\mathbf C^3$ with equation $x^3+y^3+z^3 = 0$. Then a resolution graph of $(X,0)$ consists of a single vertex with Euler class $-3$ and genus $1$. Let $L_5$ be the plumbing multilink in ${\L}_X$ whose graph is represented on figure 1. Then $L_5$ is fibred if and only if $3$ divides $n$.

\begin{figure}[!htbp]
\centering
\psfrag{-2}[c]{$-2$}
\psfrag{-3}[c]{$-3$}
\psfrag{(-1)}[c]{$(-1)$}
\psfrag{(0)}[c]{$(0)$}
\psfrag{(1)}[c]{$(1)$}
\psfrag{(2)}[c]{$(2)$}
\psfrag{(32)}[c]{$(3/2)$}
\psfrag{3}[c]{$3$}
\psfrag{L1}[c]{$L_1$}
\psfrag{L2}[c]{$L_2$}
\psfrag{L3}[c]{$L_3$}
\psfrag{L4}[c]{$L_4$}
\psfrag{L5}[c]{$L_5$}
\includegraphics[width=\textwidth]{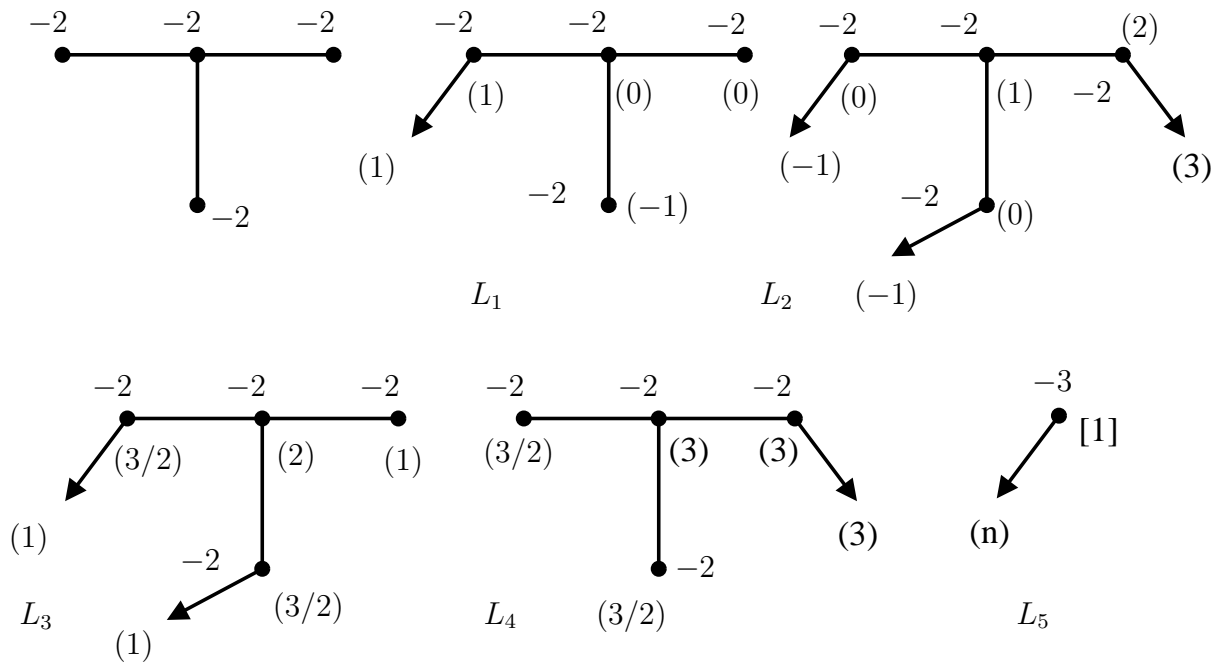}
\caption{ Multilinks}
\end{figure}

%%%%%%%%%%%%%%%%%%%%%%%%%%%%%%%%%%%%%%%%%%%%%%%%%%
%\begin{figure}[htbp]
%   \begin{center}
%      \includegraphics{fig-PS.pdf}
%   \end{center}
%   \caption{\footnotesize Figure 1}
%\end{figure}
%\begin{figure}[htbp]

%%%%%%%%%%%%%%%%%%%%%%%%%%%%%%%%%%%%%%%%%%%%%%%%%%

\begin{corollary}
\label{fibration2}
Let $(X,p)$ be a normal complex surface singularity  and let $f:(X,p) \longrightarrow
({\mathbf C},0)$ be a holomorphic germ. Let $\pi : X \to {\mathbf C}^2$ be a resolution of the holomorphic germ $fg$ such that the total transform $(fg \circ \pi)^{-1}(O)$  has normal crossings,  let $\Gamma$ be its dual graph and let $(1),\ldots,(r)$ be the  vertices of $\Gamma$. For each $i = 1,\ldots,r$, let 
$m_i^f$ (respectively $m_i^g$) be the multiplicity of $f \circ \pi$ (respectively  $g \circ \pi$) along the irreducible component of the exceptional divisor represented by $(i)$. 
Then the multilink $L_{f\overline{g}} = L_f - L_g$ is fibred if and only if for each rupture vertex $(j)$ of $\Gamma$ one has $m_j^f \neq m_j^g$.
\end{corollary}

\begin{proof}
With the notations of the theorem one has: 
$$(M_{\Gamma})^{-1}\hspace{1mm} {}^{t} b(L_f-L_g) = (M_{\Gamma})^{-1}\hspace{1mm} {}^{t} (b(L_f)-b(L_g))\,.$$

\noindent
According to ({\cite{L}, 2.6}) one has:
 $$(M_{\Gamma})^{-1}\hspace{1mm} {}^{t} b(L_f)= -{}^{t}(m^f_1,\ldots,m^f_r) \quad  \hbox{and} \quad 
 (M_{\Gamma})^{-1}\hspace{1mm} {}^{t} b(L_g)= -{}^{t}(m^g_1,\ldots,m^g_r)\,.$$
Therefore 
$$(M_{\Gamma})^{-1}\hspace{1mm} {}^{t} b(L_f-L_g) = - {}^{t} (m^f_1-m^g_1,\ldots,m^f_r-m^g_r)\,.$$ 
\qed
\end{proof}

\section{The geometry  near the multilink $L_{f \overline g}$}

In this section $(X,p)$ denotes a singularity of a normal complex surface and $f : (X,p) \to  (\mathbf{C},0) $ and $g : (X,p) \to  (\mathbf{C},0) $ 
denote two holomorphic germs without common branches. The technical point of this section is  lemma \ref{Binding}, which describes the behaviour of the real analytic germ $f \overline g : (X,p) \to  (\mathbf{C},0) $ in a neighbourhood  of its link $L_{f \overline g} \subset \L_X$.

\begin{lemma}
\label{Binding} 
Let $K$ be a component of the link $L_f \cup L_g$ and let $n$ be the multiplicity of $K$ as a component of the multilink $L_f \cup -L_g$. Then there exists a regular neighbourhood $N(K)$ 
of $K$ in $\L_X$, an orientation-preserving diffeomorphism $\tau : {\mathbf S}^1 \times {\mathbf  D}^2
\to N(K)$ such that $\tau ({\mathbf S}^1 \times \{0\}) = K$, and an integer $m \in {\mathbf Z}$ such that  for all $(t,z) \in \mathbf{S}^1 
\times ({\mathbf  D}^2 \setminus \{0\})$ we have:
$$\bigg( \frac{f \overline{g}}{|f \overline{g}
|} \circ \tau \bigg) \ (t,z)= 
\bigg( \frac{z}{\vert z \vert} \bigg)^{n} t^{m} \;.$$
\end{lemma}

\noindent
 As an immediate consequence of  \ref{Binding} and \label{fgbarre} we obtain the following :
  
\begin{corollary}
\label{CS} 
Let $(X,p)$ be a normal complex surface singularity  and let $f : (X,p) \to  (\mathbf{C},0) $ and $g : (X,p) \to  (\mathbf{C},0) $ 
be two holomorphic germs without common branches. 
If the real analytic germ $f \overline g : (X,0) \to ({\mathbf R}^2,0)$ has $0$ as an isolated critical value, then its Milnor fibration 
$\frac{f \overline g}{|f \overline g|}: \L_X \setminus 
(L_{f} \cup L_g)  \to {\mathbf S}^1$ is a fibration of the multilink $ L_f \cup -  L_g$.
\end{corollary}

\begin {remark} Notice that $L_{f \overline g}$ is  $L_f \cup L_g$ as an unoriented link, but corollary \ref{CS} states that in fact  $L_{f \overline g}$ fibres with its natural multilink structure. Thus  \ref{CS} is  much stronger than Theorem 2.1. Indeed  2.1 states that the Milnor fibration  $\frac{f \overline g}{|f \overline g|}$ is a $C^{\infty}$ locally trivial fibration of the complement of the unoriented link $(L_{f} \cup L_g) $ in $\L_X$, whereas \ref{CS} describes this Milnor fibration as the fibration of a multilink, taking into account the``multiple open-book" structure near the binding, the orientations and the multiplicities. In particular \ref{CS}  enables us to describe (using \cite{Pi1}) the genus of the Milnor fibre and most of the data determining the topology of the fibration (see for instance \cite{PS}).

\end{remark}

\begin{proof} The proof is adapted from that of (\cite{Pi2}, 3.1). Let $\pi:\tilde{X} \longrightarrow X$ be a resolution of the
holomorphic germ $fg$ such that the total
transform
$((fg) \circ \pi)^{-1}(0)$ has normal crossings. Let $S$ be a branch of the
strict transform of $f$ by $\pi$ and let $E$ be the irreducible component of the exceptional divisor 
$\pi^{-1}(p)$ which intersects $S$. We denote by $n$ the multiplicity of $f$ along the curve $S$ and by $m_f$ (respectively $m_g$) the multiplicity of $f \circ \pi$ 
(respectively $g \circ \pi$ ) along $E$. We choose local coordinates $(z_1,z_2)$ in $\tilde{X}$ such that $S \cap
E=(0,0)$, $z_2 = 0$ is a local equation for $E$ and such that in a neighbourhood of $(0,0)$ we have

$$(f \circ \pi)(z_1,z_2)=z_2^{m_f} \gamma(z_1,z_2)^n u(z_1,z_2)\,,$$
and
$$(g \circ \pi)(z_1,z_2)=z_2^{m_g} v(z_1,z_2),$$
\noindent
where $\gamma \in {\mathbf C}\{\{z_1,z_2\}\}$, $\gamma(z_1,z_2)=0$ is a
local equation of $S$ and $u$ and $v$ are units in ${\mathbf C}\{\{z_1,z_2\}\}$.
Then, locally,
$$(f\overline{g} \circ \pi)(z_1,z_2) = z_2^{m_f} \,
{\overline{z_2}}^{m_g} \, \gamma(z_1,z_2)^n \, u(z_1,z_2) \,
\overline{ v(z_1,z_2)}\,.$$

As $u$ and $v$ are units in ${\mathbf C}\{\{z_1,z_2\}\}$, one can choose two units $r$ and $s$ in ${\mathbf C}\{\{z_1,z_2\}\}$ such that 
locally, $u(z_1,z_2) = r(z_1,z_2)^n$ and  $v(z_1,z_2) = s(z_1,z_2)^n$.

Let $\alpha : {\mathbf C}^2 \longrightarrow {\mathbf C}^2$ be the map defined by
$$\alpha(z_1,z_2) \,= \, \big( \gamma(z_1,z_2)\,r(z_1,z_2) \,
\overline{ s(z_1,z_2)}\,,\,z_2\big)\,.$$

\n The jacobian determinant of $\alpha$ at the origin with respect to the
variables $(z_1,\overline{z_1},z_2,\overline{z_2})$ is:

$$|\hbox{Jac} \;\alpha (0) | = 2 i \,\Big|\frac{\partial \gamma }
{ \partial z_1}(0,0) \cdot
r(0,0) \cdot s(0,0) \Big|^2 \,;$$
but $r(0,0) \neq 0$ and $s(0,0) \neq 0$ because $r$ and $s$ are units in ${\mathbf 
C}\{\{z_1,z_2\}\}$. Moreover $\frac{\partial \gamma}
{ \partial z_1}(0,0) \neq 0$ since $S$ is transverse to $E$.
Therefore $\alpha$ is an orientation-preserving $C^{\infty}$ diffeomorphism
from an open
neighbourhood of $0_{\mathbf C^2}$ to an open neighbourhood of $0_{\mathbf C^2}$.
Replacing the local coordinates
$(z_1,z_2)$ by
$(z'_1,z_2) = \alpha (z_1,z_2)$, one obtains :

$$(f\overline{g} \circ \pi)(z'_1,z_2) = z_2^{m_f}
{\overline{z_2}}^{m_g}{z'_1}^n,$$

\noindent
where $z'_1 = 0$ is a local equation of $S$.

Now let us identify $\L_X$ with the boundary of a regular
neighbourhood $W$ of the exceptional divisor
$\pi^{-1}(p)$ in $\tilde{X}$,  defined locally by
$W=\{(z'_1,z_2) / \ / \ |z_2| \leq \eta\}$, with $\eta << 1$;  let us
study the restriction of
$\frac{f \overline{g}}  {|f \overline{g}|} \circ \pi$ to a small
tubular neighbourhood $N(K)$ of the component
$K = S
\cap \partial W$ of $L_f$, say

$$N(K) = \{ (z'_1,z_2) \  / \ |z'_1| \leq \eta', \ |z_2| = \eta \},$$

\noindent
where $\eta' \ll \eta$. Let $\tau : {\mathbf S}^1 \times {\mathbf  D}^2
\to N(K)$ be the orientation preserving diffeomorphism defined by $\tau (t,z) = (\eta' z,\eta t)$. Then for all $(t,z) \in  \mathbf{S}^1 
\times ({\mathbf  D}^2 \setminus \{0\})$ one has: 

$$\frac{f \overline{g}}{|f \overline{g}
|} \circ \pi \circ \alpha \circ  \tau \ (t,z)= \bigg( \frac{z}{\vert z \vert} \bigg)^{n} t^{m_f} \overline{t}^{m_g} =   \bigg( \frac{z}{\vert z \vert} \bigg)^{n} t^{m_f - m_g} \,.$$

\n Similar computations near a branch of the strict transform of $g$ complete
the proof.
\qed
\end{proof}

 \section{$f \overline g$ and the discriminantal ratios of $(g,f)$}

Let $(X,p)$ be a germ of  normal complex surface and let $f : (X,p) \to ({\mathbf C},0)$ and $g : (X,p) \to ({\mathbf C},0)$ be two holomorphic germs without common branches. In this section, we first give a sufficient condition for the germ $f \overline g$ to have $0$ as an isolated critical value  in terms of the discriminantal ratios of the germ $(g,f)$ with respect to the canonical complex coordinates $(u,v)$ of $(g,f)(X)$, {\it i.e.} $u=g(x)$ and $v=f(x)$ (Lemma \ref{CN}). A result due to D.T. L\^e, H. Maugendre and C. Weber (\cite{LMW}), relating the discriminantal ratios of the germ $(g,f)$ to some topological invariants of the meromorphic function $f/g$,  leads us to one of the main results of this paper, namely theorem \ref{main}.

The discriminantal ratios  are defined in \cite{LMW}  as follows. Let $C$ be the critical locus of $\pi$ and let  $\Delta = \pi(C)$ be its discriminant curve.
Let $(\Delta_{\alpha})_{\alpha \in A}$ be the set of the branches of $\Delta$ which are not the coordinates axes.

\begin{definition}  The {\it discriminantial ratio} of $\Delta_{\alpha}$ with respect to the canonical complex coordinates $u=g(x,y)$ and $v=f(x,y)$ is the  rational number   
$$d_{\alpha} = \frac{I(u=0,\Delta_{\alpha})}{I(v=0,\Delta_{\alpha})}\,,$$
where $I(-,-)$ denotes the intersection number at $0$ of complex analytic curves in ${\mathbf C}^2$.

\end{definition}

\begin{lemma}
\label{CN} 
Let $(X,p)$ be a germ of complex normal surface and let $f : (X,p) \to ({\mathbf C},0)$ and $g : (X,p) \to ({\mathbf C},0)$ be two holomorphic germs without common branches. If the discriminantal ratios of
the germ $(g,f) : (X,p) \to (\mathbf{C}^2,0)$ with respect to the complex coordinates $u=g(x)$ and $v=f(x)$ are all
different from $1$, then the real analytic germ
$f\overline{g}: (X,p) \to (\mathbf{C},0)$  has $0$ as an  isolated critical value.
\end{lemma}

\begin{proof} Let $U$ be an open neighbourhood of $p$ in $X$ such that $f$ and $g$ have no critical points in $U \setminus (fg)^{-1}(0)$. Let $x \in U \setminus (fg)^{-1}(0)$, let $\phi : V \to W$ be a parametrisation of $X \setminus \{p\}$ in a neighbourhood $W$ of $x$ in $X$, where $V$ is an open subset in ${\mathbf C}^2$ and let $y = \phi^{-1}(x)$.
As $D\phi(y) : T_y{\mathbf C}^2 \to T_xX$ is an isomorphism, then $x$ is a critical point of $f\overline{g}$ if and only if $y$ is a critical point of the map $(f\overline{g}) \circ \phi = (f \circ \phi).(\overline{g \circ \phi}) : V \to {\mathbf R}^2$.

Let us consider the two holomorphic functions  $F = f \circ \phi : V \to {\mathbf C}$ and  $G = g \circ \phi : V \to {\mathbf C}$ and let us compute the jacobian matrix of the map $(f\overline{g}) \circ \phi = F\overline{G}$ at $y = (z_1,z_2)$ with respect to the variables $z_1,\overline{z_1}, z_2,\overline{z_2}$. We decompose this map into its real and imaginary parts:
$$\big(\Re(F \overline{G}),\Im(F \overline{G})\big)=\big({1 \over
2}(F\overline{G} +  \overline{F}G),
{1 \over {2i}} (F\overline{G} -  \overline{F}G)\big) \,.$$

\noindent
Then $\,Jac(F \overline{G})(z_1,\overline{z_1},z_2,\overline{z_2})\,$ is equal to:
$$ {1 \over 2}
 \left( 
\begin{array}{cccc}
 \frac{\partial F}{\partial z_1} \overline{G} +\overline{F}  \frac{\partial G}{\partial z_1} &  
F \overline{\frac{\partial G}{\partial z_1}}  +G \overline{  \frac{\partial F}{\partial z_1}}  & 
  \frac{\partial F}{\partial z_2} \overline{G} +\overline{F}  \frac{\partial G}{\partial z_2} &  
F \overline{\frac{\partial G}{\partial z_2}}  +G \overline{  \frac{\partial F}{\partial z_2}}  \\

 -i(\frac{\partial F}{\partial z_1} \overline{G} -\overline{F}  \frac{\partial G}{\partial z_1})  &
-i(F \overline{\frac{\partial G}{\partial z_1}}  -G \overline{  \frac{\partial F}{\partial z_1}} ) & 
  -i(\frac{\partial F}{\partial z_2} \overline{G} -\overline{F}  \frac{\partial G}{\partial z_2})  &
-i(f \overline{\frac{\partial G}{\partial z_2}}  -G \overline{  \frac{\partial F}{\partial z_2}} ) \\
\end{array}
\right) \;.
$$

\n Therefore the rank of
$Jac(F
\overline{G})(z_1,\overline{z_1},z_2,\overline{z_2})$  at $(z_1,z_2) \in V$ is not
maximal if and only if the following four conditions hold :

$$(1) \hspace{.5in}     FG \big( {\frac{\partial F}{\partial z_1}} \frac{\partial G}{\partial
z_2}  -  \frac{\partial F}{\partial z_2} \frac{\partial G}{\partial z_1}  \big) 
=  0 \;;
$$

$$(2) \hspace{.5in} \mid G \frac{\partial F}{\partial z_1}\mid   = \mid F  \frac{\partial G}{\partial z_1}  \mid \;;$$

$$(3) \hspace{.5in}  \mid G \frac{\partial F}{\partial z_2}\mid  =  \mid F  \frac{\partial G}{\partial z_2}\mid \;; $$

$$(4) \hspace{.5in} \mid G \mid^2  \frac{\partial F}{\partial z_1} \overline{\frac{\partial F}{\partial z_2}}  = 
\mid F \mid^2  \frac{\partial G}{\partial z_1} \overline{\frac{\partial G}{\partial z_2}}\;. $$

Assume that $(z_1,z_2) \in V \setminus (FG)^{-1}(0)$ is such a point.  Then (1) implies

$$
 {\frac{\partial F}{\partial z_1}} (z_1,z_2)  \frac{\partial G}{\partial
y} (z_1,z_2) -  \frac{\partial F}{\partial z_2}(z_1,z_2)  \frac{\partial G}{\partial z_1}(z_1,z_2) =0 \;,
$$

\noindent
{\it i.e.} $(z_1,z_2)$ belongs to the critical 
locus of the germ $(G,F)$, and then $\phi(z_1,z_2)$ belongs to the critical locus $C$ of $(G,F)$.

Now assume that $f\overline g$ does not have $0$ as an isolated critical value. 
Then there exists a branch $\gamma$ of
$C$ in $U $ which is not a branch of ${fg}^{-1}(0)$, and a sequence
$\{x_n\}_{n \in \mathbf{ N}}$ of points of $\gamma \setminus \{p\}$ converging
to $p$ which satisfy that  for each $n \in {\mathbf N}$ and for every parametrisation $\phi : V \to W$ of $X \setminus \{p\}$ such that $x_n \in W$, the rank of
$Jac(f \overline{g} \circ \phi )(\phi^{-1}(x_n))$ is not maximal. Let us fix a  parametrisation $\phi : V \to W$
of $X\setminus \{p\}$ on $W$ such that $W$ is a contractible open subset in $X \setminus \{p\} $ such that $\gamma \setminus \{p\} \subset W$. Let us set $\beta = \phi^{-1}(\gamma \setminus \{p\})$,  $F = f \circ \phi$, $G =  g \circ \phi$, and for each $n \in {\mathbf N}, (z_{n,1},z_{n,2}) = \phi^{-1}(x_n)$ and $(u_n,v_n)=(g(x_n),f(x_n) $).

Let $\alpha \in A$ be such that $\Delta_{\alpha} =
(g,f)(\gamma)$. According to the definition of the discriminantal ratio $d_{\alpha}$,  the curve $\Delta_{\alpha}$ admits a Puiseux expansion of the form :
$$u= \Psi(v) =  v^{d_{\alpha} }\bigg(a+\sum_{k \in {\mathbf N}^{\ast}}
b_k v^{k \over m }\bigg)\;.$$

Then $G_{\mid \beta} = (\Psi \circ F
)_{\mid \beta}$, and
for each point $(z_1,z_2) $ of the curve $ \beta $, the linear maps
$DG(z_1,z_2)$ and $D(\Psi \circ F )(z_1,z_2)$ coincide
on the complex tangent line
${ T_{(z_1,z_2)}\beta} $.  Set $v=F(z_1,z_2)$ and let
$h \in { T_{(z_1,z_2)}\beta} \setminus \{0\}$. Then $DG(z_1,z_2)\cdot h =\Psi'(v)
DF(z_1,z_2) \cdot h$. As $(z_1,z_2)$ belongs to the jacobian
locus of $(G,F)$, for all $ (z_1,z_2) \in \beta$ one  obtains:

$$ 
 \bigg(  \frac{\partial G}{\partial z_1}(z_1,z_2),\frac{\partial G}{\partial z_2 }(z_1,z_2) \bigg) =
\Psi'(v) \bigg(  \frac{\partial F}{\partial z_1}(z_1,z_2),\frac{\partial G}{\partial z_2}(z_1,z_2) \bigg)\;.
$$

\n This holds in particular for each$(z_{n,1},z_{n,2}) , n \in {\mathbf N}$.

Since $f$ and $g$ have no singular point in  $  U \setminus (fg)^{-1}(0)$, then for each $n \in {\mathbf N}$, either $ \big(  \frac{\partial F}{\partial z_1}(z_{n,1},z_{n,2}), \frac{\partial G}{\partial z_1}(z_{n,1},z_{n,2}) \big)
\neq (0,0)$
 or  $ \big(  \frac{\partial F}{\partial z_2}(z_{n,1},z_{n,2}), \frac{\partial G}{\partial z_2}(z_{n,1},z_{n,2}) \big)
\neq (0,0)$
.
 Then, perhaps after replacing $\{x_n\}_{n \in {\mathbf N}}$ by a subsequence, one can assume that
either 

\v
(a)  for all $n \in \mathbf{ N}$ one has   $\big(  \frac{\partial F}{\partial z_1}(z_{n,1},z_{n,2}), \frac{\partial G}{\partial z_1}(z_{n,1},z_{n,2}) \big)
\neq (0,0)$;

\noindent
or 

\v 
(b)  for all $n \in \mathbf{ N}$ one has  $ \big(  \frac{\partial F}{\partial z_2}(z_{n,1},z_{n,2}), \frac{\partial G}{\partial z_2}(z_{n,1},z_{n,2}) \big)
\neq (0,0)\;.$ 
 
\v
If (a) holds then (2) implies that for all $n \in \mathbf{ N}$ one has $ \ |\Psi(v_n)| = |v_n \Psi'(v_n)|$,  where $v_n = f(x_n)$, {\it i.e.}

$$ \ |a d_{\alpha}+ \sum_{k \in
{\mathbf{ N}}^{\ast}} b_k(d_{\alpha} + \frac{k}{ m}){v_n}^{\frac{k}{m}}|=|a+ \sum_{k \in {\mathbf N}^{\ast}} 
b_k {v_n}^{\frac{k}{m}}|\,,$$
for all $n \in \mathbf{ N}$. 
Taking the limit when $n \to \infty$ leads  to $|a d_{\alpha}|=|a|$, {\it i.e.} 
    $d_{\alpha}=1 $.
                                  
Condition (b) also leads to  $d_{\alpha}=1 $.
                 \qed
\end{proof}

 The following result is a generalization of (\cite{Pi2}, theorem 5.1), which treated the particular case where $X={\mathbf C}^2$ and  $f$ and $g$ have an isolated singularity at $0$. 
  
\begin{theorem} 
\label{main}
Let $(X,p)$ be a singularity of normal surface and let $f: (X,p) \to ({\mathbf C},0)$ and $g: (X,p) \to ({\mathbf C},0)$ be two holomorphic germs without common branches. Then the following conditions are equivalent
\begin{description}
\item[( i)] The real analytic germ $f \overline g : (X,p) \to (\mathbf R^2,0)$ has $0$ as an isolated critical value.
\item[( ii)] The multilink $L_f -L_g$ is fibred. 
\item[( iii)] Let $\pi : \tilde{X} \to X$ be any resolution of the holomorphic germ $fg: (X,p) \to ({\mathbf C},0)$. For each rupture vertex $(j)$ of the dual graph of $\pi$, $m_j^f - m_j^g \neq 0$.
 
\end{description}
Moreover, if these conditions hold, then the Milnor fibration $\frac{f \overline{g} }{ |{f \overline {g}}|} : 
\L_X \setminus (L_f \cup L_g)$ of ${f \overline g}$ is a fibration of the multilink $L_f -L_g$.
\end{theorem}

Notice that this is Theorem 2 in the introduction. 

The key point of the proof of (\cite{Pi2}, 5.1)   is the theorem 1.1 of \cite{Ma}, which relates the determinantal ratios of the germ $(g,f) : ({\mathbf C}^2,0) \to ({\mathbf C}^2,0)$ with some topological invariants of the meromorphic function $f/g$. The key point of the proof of  theorem \ref{main},  is the following generalization of (\cite{Ma}, theorem 1.1).
 
 \begin{definition} Let $\pi : \tilde{X} \to X$ be any resolution of the holomorphic germ $fg$, let $\Gamma$ be the associated resolution graph, and let $(1),\ldots,(r)$ be the vertices of $\Gamma$. For each $i = 1,\ldots,r$, let 
$m_i^f$ (respectively   $m_i^g$) be the multiplicity of $f \circ \pi$ (respectively   $g \circ \pi$) along the irreducible component of the exceptional divisor $\pi^{-1}(p)$ represented by $(i)$. The {\it contact quotient} associated with $(i)$ is the rational number $\frac{m_i^f}{m_i^g}$.
\end{definition}

\begin{theorem} (\cite{LMW}, Theorem 0.3) 
\label{LMW}
The set of the discriminantal ratios of the germ $(g,f) : (X,p) \to ({\mathbf C}^2,0)$ equals the set $T$ of  contact quotients of $(g,f)$ associated to the rupture vertices of a resolution graph of the germ $fg$. 
\end{theorem}
 \noindent
Notice that the set $T$ does not depend on the choice of the resolution.

\begin{proof}  {\bf of \ref{main}} By   \ref{CS}, if $f \overline g $ has $0$ as an isolated critical value, then the multilink $L_f -L_g$ is fibred.
  
Let $\pi : \tilde{X} \to X$ be a resolution of the holomorphic germ $fg$ and let $\Gamma$ be its dual graph. Let $(1),\ldots,(r)$ be the  vertices of $\Gamma$.  Let us denote by $\mathcal R$ the set of rupture vertices of $\Gamma$. 

Assume that the plumbing multilink $L_f -L_g$ is fibred. Then, according to  \ref{fibration2}, for each vertex $(j) \in \mathcal{R}, m_j^f \neq m_j^g$, {\it i.e.} $\frac{ m_j^f}{m_j^g} \neq 1$. But $\frac{ m_j^f}{m_j^g}$ is the contact quotient of $(g,f)$ associated with the vertex $(j)$. Then, according to \ref{LMW},  the set of the discriminantal ratios of the germ $(g,f)$ does not contain $1$. Therefore, by  lemma \ref{CN}, $f \overline g $ has $0$ as an isolated critical value.
\qed
\end{proof}

\section{Realization of fibred multilinks by singularities : results and open questions}

A natural family of questions, which can be called ``realization questions", arises in the study of the topological properties of singularities. Roughly speaking, a realization question asks wether a given topological object (e.g. a link in $\mathbf S^3$, a fibred multilink, a diffeomorphism of a surface, ...) is realized by a singularity (e.g. as the link of an analytic germ of function, by a Milnor fibration, as the monodromy of a Milnor fibration, ...).

In \cite{AK}, it is proved that all knots (in fact all links) are algebraic in the following sense. Let $U\subset S^{k-1}$ be a compact smooth submanifold (with boundary) of codimension $\geq 1$ with trivial normal bundle. There is a real algebraic set 
$Z \subset R^k$ with an isolated singularity at the origin such that for all small
 $\varepsilon >0$ the pair $(S^{k-1}_\varepsilon ,S^{k-1}_\varepsilon  \cap Z)$ 
is diffeomorphic to $(S^{k-1},\partial U)$. If $U$ is a Seifert surface of the knot $K$ in $S^3$ the theorem yields the result stated as title of the paper, i.e. ``All knots are algebraic'' (see also \cite{MVS}).
More precisely, for each link $L$ in $\mathbf S^3$, the algebraic set is the  fibre $f^{-1}(0)$ of a  polynomial function $f:{\mathbf R}^4 \longrightarrow {\mathbf R}^2$ such that $f(0)=0$, and  the pair $(S^3,L)$ is diffeomorphic to the pair $({\mathbf S}_{\epsilon}, L_f)$ where $L_f = f^{-1}(0) \cap {\mathbf S}_{\epsilon} , \ 0 < \epsilon <<1$ is the link associated with the analytic germ  $f:({\mathbf R}^4 ,0) \longrightarrow ({\mathbf R}^2,0)$.

This section is devoted to the following realization question : given a plumbing $3$-manifold $M$ and a fibred multilink $L=n_1 k_1 \cup \ldots \cup n_l  K_l$ in $M$, does there exist some   holomorphic or real analytic germ $f:(X,p) \to ({\mathbf C},0)$ from a germ $(X,p)$ of a complex normal surface $(X,p)$ such that $(M,L)$ is diffeomorphic to $({\cal L}_X,L_f)$ ?

We recall the convention introduced in Section 4 : in this paper,  the Seifert pieces of any  plumbing manifold have all orientable basis and the plumbing operations are all positive. 

In \cite{Pi3}, using a realization result of G. Winters \cite{Wi} for degenerating families of complex, we prove the following realization result for positive multilinks : 

\begin{theorem} \label{real1}
Let $M$ be a plumbing  manifold and let $L=n_1 k_1 \cup \ldots \cup n_l  K_l$ be a multilink in $M$ such that $ n_i >0$ for all $i = 1,\ldots,l$.
If $L$ is fibred, then there exists a  normal complex surface singularity  $(X,p)$ and a holomorphic germ $f : (X,p) \to ({\mathbf C},0)$ such that the pair $(M,L)$ is diffeomorphic to the pair $({\cal L}_X,L_f) $.
\end{theorem}

\begin{remark}
The condition ``L is fibred" is easy to check from any plumbing graph $\Gamma$ of $(M,L)$ using Theorem \ref{fibration1}. Also notice that the conditions 1) and 2) of \ref{fibration1} imply that the intersection matrix $M_{\Gamma}$ is negative definite, {i.e.} Grauert's condition for $M$ to be homeomorphic to the link of a normal surface singularity.
\end{remark}

As an immediate consequence of Theorem \ref{real1} and Theorem \ref{main} one obtains the following realization result, which is Theorem 3 in the introduction :

\begin{theorem} 

 Let $M$ be a plumbing  $3$-manifold  and let $L_1$ and $L_2$ be two plumbing multilinks in $M$ with positive multiplicities which verify the following conditions : 
 
1.  the three multilinks $L_1 \cup (- L_2)$, $L_1$ and $L_2$ are fibred;

2. $M$ is either $\S^3$ or 
 homemorphic to the link ${\cal L}_X$ of a  taut  surface singularity $(X,p)$.
 \n 
Then, there exist two holomorphic germs $f,g:(X,p) \to ({\mathbf C},0)$ without common branches such that $L_1$ (and $L_2$) is the multilink of $f$ (respectively $g$); $L_1 \cup (-L_2)$ is the multilink associated to the real analytic germ $f \overline{g} : (X,p) \to ({\mathbf C},0)$; and
$$\frac{f \overline g}{\vert f \overline g \vert} \,:\, 
{\cal L}_X \setminus (L_1 \cup (- L_2)) \longrightarrow \S^1\;,$$
is a fibre bundle that realizes $L_1 \cup (- L_2)$ as a fibred multilink.
\end{theorem}

Notice that if $L_1 \cup (- L_2)$ and either  $L_1$ or $L_2$ are fibred, then the third link is also fibred. 
As already mentioned, the condition 1 implies that $M$ is homeomorphic to the link ${\cal L}_X$ of a singularity of a normal complex surface, but in general $(X,p)$ is not taut, {\it i.e.} the analytical type of the germ $(X,p)$ is not unique.

\v
In the following two situations  the realization problem remains open : 

\v
\noindent
1) {\it The condition 1 in 7.2 
 is verified and  $M$ is homeomorphic to the link ${\cal L}_X$ of a singularity of a normal complex surface which is not taut}. In this case, Theorem \ref{real1} implies that the multilinks 
$L_1$ and $L_2$ are  realized as the multilinks of two holomorphic germs $f : (X_1,p_1) \to  ({\mathbf C},0)$ and  $g : (X_2,p_2) \to  ({\mathbf C},0)$ respectively, but in general, $X_1$ and $X_2$ cannot be chosen with the same analytic type. Therefore the germ  $f \overline  g$ may not be defined and, {\it a priori}, 
the multilink link 
$L_1 \cup -L_2$ might not be realizable as the multilink of a real analytic germ $f \overline  g$ by this method.

\v
\noindent
2) {\it The link $L_1 \cup (- L_2)$ is fibred, but the multilinks $L_1$ and $L_2$ are not fibred} (e.g.  the multilink $L_3 \cup (-L_4)$ of example 4.9).  According to \cite{Pi3}, since $L_i$ has positive multiplicities, then $m_j^{(i)} >0$
for all $j \in \{1,\ldots,r\}$,  where 
$$ {}^{t}(m_1^{(i)},\ldots,m_r^{(i)},) =-(M_{\Gamma})^{-1}\hspace{1mm} {}^{t} b(L_i).$$
Therefore $L_i$ is not fibred if and only if there exists $j \in \{1,\ldots , r \}$ such that $m_j^{(i)} \not\in {\mathbf Z }$. Since $L_1 \cup -L_2$ is fibred by hypothesis,
 then for each such $j$, $m_j^{(1)} - m_j^{(2)} \in {\mathbf Z }$.

%---------------------

\end{document}